\documentclass[11pt]{article}
\usepackage{epsfig}

\def\be{\begin{equation}}
\def\ee{\end{equation}}
\def\bea{\begin{eqnarray}}
\def\eea{\end{eqnarray}}
\def\bes{\begin{eqnarray*}}
\def\ees{\end{eqnarray*}}

\def\nn{\nonumber}
\def\<{\langle}
\def\>{\rangle}
\def\lb{\label}
\def\bs{\setminus}
\def\pt{\partial}

\def\R{{\bf R}}

\def\Z{{\bf Z}}
\def\N{{\bf N}}

\def\Q{{\bf Q}}

\def\aa{{\alpha}}

\def\ga{{\gamma}}
\def\Ga{{\Gamma}}
\def\ka{{\kappa}}
\def\th{{\theta}}

\def\ep{{\epsilon}}
\def\lm{{\lambda}}
\def\Lm{{\Lambda}}

\def\sg{{\sigma}}

\def\ev{{\it ev}}

\def\ker{{\rm ker}}

\def\index{{\it index}}
\def\nullity{{\it nullity}}

\def\Sp{{\rm Sp}}
\def\mod{{\rm mod}}

\def\ol#1{\overline{#1}}  
\def\td#1{\tilde{#1}}
\def\hb{\vrule height0.18cm width0.14cm $\,$}
\def\mapright#1{\smash{\mathop{\longrightarrow}\limits^{#1}}}

\title{The existence of two closed geodesics on\\
every Finsler 2-sphere}
\author{Victor Bangert\thanks{Partially supported by DFG-Forschergruppe
`Nonlinear Partial Differential Equations: Theoretical and Numerical
Analysis'.} \quad and \quad
Yiming Long\thanks{Partially supported by the 973 Program of MOST,
NSFC, Cheung Kong Scholarship, MCME, SRFDP of MOE of China, Li Ka Shing
Foundation, Shiing-Shen Chern Foundation, LPMC, and Nankai University. }}
\date{}

\begin{document}

\maketitle

\begin{abstract}
{\it In this paper, we prove that for every Finsler metric on $S^2$ there exist
at least two distinct prime closed geodesics. For the case of the two-sphere,
this solves an open problem posed by D. V. Anosov in 1974. }
\end{abstract}



\renewcommand{\theequation}{\thesection.\arabic{equation}}
\renewcommand{\thefigure}{\thesection.\arabic{figure}}

\setcounter{equation}{0}
\section{Introduction}

The study of closed geodesics on spheres is a classical and
important problem both in dynamical systems and differential
geometry. The results of J. Franks \cite{Fra1} in 1992 and V.
Bangert \cite{Ban3} in 1993 prove that for every Riemannian metric
on $S^2$ there exist infinitely many geometrically distinct closed
geodesics. In 1973, remarkable irreversible Finsler metrics on
$S^2$ were constructed by A. Katok in \cite{Kat1} which possess
precisely two distinct prime closed geodesics (cf. \cite{Zil2} for
further explanation), which are inverse curves of each other.

Based on this result, D. V. Anosov in the I.C.M. of 1974 (cf. \cite{Ano1})
wrote: "For the $n$-dimensional sphere $S^n$, Katok's example gives an
irreversible Finsler metric, arbitrarily near to the 'standard' metric
(to the metric of constant curvature) which has $2[(n-1)/2]$ closed
geodesics. This number coincides with the lower bound which one naturally
expects for irreversible Finsler metrics on $S^n$ and which can be proved
for metrics sufficiently near to the 'standard' metric' ", and
mentioned that the "basic difficulty" is due to the multiple
covering phenomenon. Note that according to \cite{Kat1} and \cite{Zil2}
with the usual interpretation of $[\,\cdot\,]$, the number of distinct
closed geodesics on Katok's examples of Finsler metrics on $S^n$ is
actually $2[(n+1)/2]$. Note also that in \cite{Zil2} of 1982, W. Ziller
asks if $n$ is a lower bound for the number of distinct closed geodesics
on any Finsler sphere $S^n$ . In this paper we resolve these problems
for Finsler metrics on $S^2$. 

{\bf Theorem 1.1.} {\it For every Finsler metric on the 2-sphere,
there exist at least two closed geodesics. }

Note that in Theorem 1.1 the Finsler metric is not assumed to be reversible.
In particular, if $c:S^1=\R / \Z \rightarrow S^2$ is a closed geodesic then its inverse curve
$c^{-1}$, defined by $c^{-1}(t)=c(1-t)$, will not be a geodesic in general.
If it is, it is counted as a second closed geodesic, as in Katok's examples.

For the definition of Finsler metrics and their geodesics we refer to \cite{BCS1} and \cite{She1}.

According to the classical theorem of Lyusternik-Fet \cite{LyF1}
from 1951, there exists at least one closed geodesic on every
compact Riemannian manifold. The proof of this theorem, see
e.g.~\cite{Ban2} or \cite{Kli2}, is variational and carries over
to the Finsler case. The authors are aware of only few results on
the existence of more than one closed geodesic on Finsler
2-spheres. In \cite{Rad1} of 1989, H.-B. Rademacher proved the
existence of a second closed geodesic on a Finsler 2-sphere
 provided all its closed geodesics - including iterates - are non-degenerate.
If, in addition, for every hyperbolic closed geodesic the stable
and unstable manifolds intersect transversely, then the results by
H. Hofer, K. Wysocki and E. Zehnder in \cite{HWZ2} from 2003 imply
the existence of either two or infinitely many closed geodesics.
This alternative also holds for Finsler metrics with $K\ge 1$ for
which every geodesic loop is longer than $\pi$. This was proved in
2006 by A. Harris and G. P. Paternain \cite{HaP1}. Their proof is
based on \cite{HWZ1}. Recently, H.-B. Rademacher \cite{Rad3} proved
the existence of two closed geodesics on Finsler 2-spheres satisfying
a pinching condition on the flag curvature.

It is tempting to try to prove Theorem 1.1 along the lines of the
proof of the celebrated theorem of Lyusternik-Schnirelmann on the
existence of three closed geodesics without self-intersections on
every Riemannian 2-sphere, cf. \cite{LyS1}, \cite{Lyu1},
\cite{Bal1}, \cite{Jos1}, \cite{Gra1} and \cite{Tai1}. This would
be possible if one could find an energy-decreasing deformation on
the space of closed curves without self-intersections that works
for irreversible Finsler metrics on $S^2$. As observed by H.-B.
Rademacher, see p.82 of \cite{Rad1}, Katok's examples show that
such a deformation does not exist. Instead, our proof is by
contradiction. We assume that there is only one closed geodesic
$c$ and make a case by case study according to the different
symplectic normal forms of the linearized Poincar\'{e} map of $c$.
For more details see the outline of the proof in the next section.

Using Legendre transformation one can reformulate Theorem 1.1 as a result on
convex Hamiltonian systems on the cotangent bundle $T^{*}S^2$ of $S^2$ as follows.

{\bf Theorem 1.2.} {\it Let $H:T^*S^2\rightarrow\R$ be smooth and assume that the
restrictions of $H$ to the fibers of the cotangent bundle $T^*S^2\rightarrow S^2$
have positive definite Hessian everywhere and attain their minima. Let $r$ be a
real number such that the sublevel $H^{-1}((-\infty, r))\subseteq T^*S^2$ contains
the zero section. Then the Hamiltonian system $X_H$ determined by $H$ and the
standard symplectic structure on $T^*S^2$ has at least two periodic orbits on
the level surface $H^{-1}(r)$.}

\setcounter{equation}{0}
\section{Outline of the proof}

We first explain how we count closed geodesics on a Finsler
manifold $(M, F)$. If $c:S^1=\R/ \Z\rightarrow M$ is a closed
geodesic, then so are its {\it iterates} $c^m:S^1\rightarrow M$
defined by $c^m(s)=c(ms)$, for all positive integers $m$. The
closed geodesic $c$ is called {\it prime} if there does not exist
$\tilde{c}:S^1\rightarrow M$ such that $c=\tilde{c}^m$ for some
$m\ge 2$. Two prime closed geodesics $c$ and $d$ are {\it distinct}
if they do not only differ by translation of parameter,
i.e., if there does not exist $\theta\in\R / \Z$ such that
$c(t)=d(t+\theta)$ for all $t\in\R / \Z$. Now the number of closed
geodesics of $(M,F)$ is defined as the (possibly infinite) number
of distinct prime closed geodesics of $(M,F)$. In order to prove
Theorem 1.1 we argue by contradiction and assume that there is a
Finsler 2-sphere $(S^2, F)$ with only one prime closed geodesic
$c$. The closed geodesics $\{c^m|m\ge 1\}$ are critical points of
the Finsler energy functional $E:\Lambda\rightarrow\R$,
$E(\gamma)=\frac{1}{2}\int\nolimits_{S^1}F^2(\dot{\gamma}(t))dt$
on the space $\Lambda=\Lambda S^2$ of closed $H^1$-curves. As
critical points of $E$ the closed geodesics have an index
$$ i(c^m)=\index(D^2 E(c^m)) $$
and a nullity
$$ \nu(c^m)=\nullity (D^2 E(c^m))-1. $$
As proved in \cite{Lon1} (cf. \cite{Lon3}), there are nine possible cases for
the sequences $\{i(c^m)\}_{m\ge 1}$ and $\{\nu(c^m)\}_{m\ge 1}$ depending on
the different symplectic normal forms of the linearized Poincar\'{e} map of
$c$. Note that here by \cite{LLo2} the iteration formulae in \cite{Lon1} work
for Morse indices of closed geodesics on Riemannian and Finsler manifolds.
In most of these nine cases we show that Morse theory necessitates the
existence of an additional closed geodesic besides $c$. Here, the two
non-degenerate cases, in which $\nu(c^m)\equiv 0$, have already
been treated in \cite{Rad1}. In some of the degenerate cases the
results by N. Hingston, \cite{Hin1} and \cite{Hin2}, are of great
use. In others, we can use arguments of Gromoll-Meyer type, cf.
\cite{GrM1} and \cite{GrM2}, to obtain contradictions to the Morse
inequalities with $\Q$-coefficients. Here we need very precise
information on the homology created by iterated closed geodesics in the
loop space. This is the content of Sect. 3. It relies on H.-B. Rademacher's
Habilitationsschrift \cite{Rad2}, while part of it is new.
In some of these cases the existence of a
second closed geodesic also follows from the recent paper \cite{Rad3}
by H.-B. Rademacher. There remains one subcase, treated in
Sect. 10, in which the preceding methods fail. To obtain a
contradiction in this subcase, we use on the one hand a result by
H.-B. Rademacher from \cite{Rad2} on the mean index of $c$,
cf. Theorem 5.2 below, and on the other hand a detailed study on how
the first and second homologies of the free loop space modulo point
curves should be created by the closed geodesic and its iterates,
cf. Sect. 10.


Contents
\begin{itemize}
\item[1.] Introduction
\item[2.] Outline of the proof
\item[3.] Critical modules of iterated closed geodesics
\item[4.] Classification of closed geodesics on $S^2$
\item[5.] Rationally elliptic degenerate saddle closed geodesics on $S^2$
\item[6.] Two theorems by N.~Hingston
\item[7.] Homology of $(\Lambda S^2, \Lambda^0 S^2)$ and first consequences
\item[8.] Cases with eigenvalue $1$
\item[9.] Cases with eigenvalue $-1$
\item[10.] Case CG-7 of a rationally elliptic closed geodesic
\end{itemize}

We add some remarks concerning notations in this paper.

For $a\in\R$ we set $[a]=\max\{k\in\Z|k\le a\}$.
By $\N$ we denote the set of positive integers.
Throughout the paper, homology modules will be with respect to coefficients in $\Q$.
This allows us to use the transfer theorem for the homology of spaces with
$\Z_m$-actions, see Lemma 3.6.

\setcounter{equation}{0}
\section{Critical modules of iterated closed geodesics}

In this section we review part of the Morse theory of the energy functional
on the free loop space of a Finsler manifold. Some facts that we need do not
yet exist in the literature, and so we present some of the details. On a
compact Finsler manifold $(M,F)$ we choose an auxiliary Riemannian metric.
This endows the space $\Lambda=\Lambda M$ of $H^1$-maps
$\gamma:S^1\rightarrow M$ with a natural structure of Riemannian Hilbert
manifold on which the group $S^1=\R/\Z$ acts continuously by isometries. Here
a map $c:S^1\to M$ is $H^1$ if it is absolutely continuous and the derivative
$\dot{c}(t)$ is square integrable, cf. \cite{Kli2}, Chapters 1 and 2. The $S^1$-action is defined
by translating the parameter, i.e.
$$  (s\cdot\gamma)(t)=\gamma(t+s) $$
for all $\gamma\in\Lambda$ and $s,t\in S^1$. The Finsler metric $F$ defines an
energy functional $E$ and a length functional $L$ on $\Lambda$ by
\be
E(\gamma)=\frac{1}{2}\int_{S^1}F(\gamma(t),\dot{\gamma}(t))^2 \mbox{dt}, \quad
L(\gamma) = \int_{S^1}F(\gamma(t), \dot{\gamma}(t)) \mbox{dt}.  \lb{2.1}
\ee
Both functionals are invariant under the $S^1$-action.
For $\kappa\in\R$ we set $\Lambda^{\kappa}=\{\gamma\in\Lambda|E(\gamma)\le\kappa\}$ and
$\Lambda^{\kappa_-} = \{\gamma\in\Lambda|E(\gamma)<\kappa\}$.

The critical points of
$E$ of positive energy are precisely the closed geodesics $c:S^1\rightarrow M$
of the Finsler structure. If $F$ is not Riemannian, then due to the
non-differentiability of $F^2$ on the zero section, the energy
$E:\Lambda\rightarrow[0,\infty)$ is not smooth, but only of class
$C^{1,1}$, cf. \cite{Mer1}. If $c\in\Lambda$ is a closed geodesic, then $c$
is a regular curve, i.e. $\dot{c}(t)\not= 0$ for all $t\in S^1$, and this
implies that the second differential $E''(c)$ of $E$ at $c$ exists.

As usual we define the index $i(c)$ of $c$ as the maximal dimension of
subspaces of $T_c \Lambda$ on which $E''(c)$ is negative definite, and
the nullity $\nu(c)$ of $c$ so that $\nu(c)+1$ is the dimension of the
null space of $E''(c)$. The relations between $E''(c)$, the index form and Jacobi
fields are analogous to the Riemannian case, see e.g.~\cite{She1}.

For $m\in\N$ we denote the $m$-fold iteration map
$\phi_m:\Lambda\rightarrow\Lambda$ by
\be \phi_m(\ga)(t)=\ga(mt) \qquad \forall\,\ga\in\Lm, t\in S^1.
\lb{2.3}\ee
We also use the notation $\phi_m(\gamma)=\gamma^m$. Note that $\phi_m$
is an embedding satisfying
\be  E\circ\phi_m=m^2 E  \lb{2.4}\ee
and
\be  \phi_m(s\cdot\gamma)=\frac{s}{m}\cdot\phi_m(\gamma) \qquad
            \forall \,\gamma\in\Lambda, s\in S^1. \lb{2.5}\ee
According to (\ref{2.4}), if $c$ is a closed geodesic and
$\xi$, $\eta\in T_c\Lm$, then
\be  E''(c^m)(D\phi_m(\xi), D\phi_m(\eta))=m^2 E''(c)(\xi,\eta).
      \lb{2.6}\ee
Since the null space of $E''(c)$ is the space of periodic Jacobi fields
along $c$, one easily concludes:

{\bf Lemma 3.1.} {\it
$D\phi_m(c)$ maps the null space of $E''(c)$ injectively into the null space
of $E''(c^m)$. }

If $\gamma\in\Lambda$ is not constant then the multiplicity $m(\gamma)$ of
$\gamma$ is the order of the isotropy group
$\{s\in S^1\mid s\cdot\gamma=\gamma\}$. If $m(\gamma)=1$, then $\gamma$ is
called prime. Hence $m(\gamma)=m$ if and only if there exists a prime
curve $\tilde{\gamma}\in\Lambda$ such that $\gamma=\tilde{\gamma}^m$.

For a closed geodesic $c$ we set
$$ \Lm(c)=\{\ga\in\Lm\mid E(\ga)<E(c)\}. $$
If $A\subseteq\Lm$ is invariant under some subgroup $\Gamma$ of $S^1$,
we denote by $A/\Gamma$ the quotient space of $A$ with respect to the
action of $\Gamma$.

Using singular homology with rational coefficients we will consider
the following critical $\Q$-modules of a closed geodesic $c\in\Lm$:
\be \left\{\matrix{
C_*(E,c) &=& H_*\left(\Lambda(c)\cup\{c\}, \Lambda(c)\right), \cr
C_*(E,S^1\cdot c) &=& H_*\left(\Lambda(c)\cup S^1\cdot c, \Lambda(c)\right), \cr
\ol{C}_*(E,c)
&=& H_*\left((\Lm(c)\cup S^1\cdot c)/S^1, \Lm(c)/S^1\right). \cr}\right.
        \lb{2.8}\ee

In order to relate the critical modules to the index and nullity of $c$ we
would like to use the results by D. Gromoll and W. Meyer from \cite{GrM1},
\cite{GrM2}. Unfortunately, this is not directly possible, since in general
the functional $E$ will not be of class $C^2$ in any neighborhood of $c$.
Following \cite{Rad2}, Section 6.2, we will evade this problem by
introducing finite-dimensional approximations to $\Lambda$. We choose an
arbitrary energy value $a>0$ and $k\in\N$ such that every $F$-geodesic of
length $<\sqrt{2a/k}$ is minimal. Then
$$ \Lm(k,a)=\left\{\ga\in\Lm\mid E(\ga)<a \mbox{ and }
  \ga|_{[i/k,(i+1)/k]}\mbox{ is an $F$-geodesic for }i=0,\ldots,k-1\right\} $$
is a $(k\cdot\dim M)$-dimensional submanifold of $\Lambda$ consisting of
closed geodesic polygons with $k$ vertices. The set $\Lambda(k,a)$ is
invariant under the subgroup $\Z_k$ of $S^1$. The important point is that the
energy functional $E$ is smooth on the open and dense subset of $\Lambda(k,a)$
that consists of all polygons for which no two consecutive vertices coincide.
Moreover the closed geodesics in
$\Lambda^{a_-}=\{\gamma\in\Lambda\mid E(\gamma)<a\}$ are precisely the
critical points of $E|_{\Lm(k,a)}$, and for every closed geodesic
$c\in\Lm(k,a)$ the index of $(E|_{\Lm(k,a)})''(c)$ equals $i(c)$
and the null space of $(E|_{\Lm(k,a)})''(c)$ coincides with the
nullspace of $E''(c)$, cf. \cite{Rad2}, p.51. Finally, there exists a
$\Z_k$-equivariant, energy non-increasing deformation retraction
\be r:\Lambda^{a_-}\times[0,1]\rightarrow\Lambda^{a_-} \lb{2.9}\ee
of $\Lambda^{a_-}$ to $\Lambda(k,a)$, cf.~\cite{Rad2}, Section 6.2.
It is defined by
$$ \begin{array}{rcl}
r(\gamma, u)|_{[i/k, (i+u)/k]} &=& \mbox{ the minimal geodesic from }
\gamma(i/k) \mbox{ to } \gamma((i+u)/k), \\
r(\gamma, u)|_{[(i+u)/k, (i+1)/k]} &=& \gamma|_{[(i+u)/k, (i+1)/k]},
\end{array}
$$
for $\gamma\in\Lambda^{a_-}$, $u\in[0,1]$ and $i=0,\ldots,k-1$.

Throughout the paper we will assume that each closed geodesic
$c\in\Lambda$ satisfies the following isolation condition:

{\bf (Iso) For all $m\ge 1$ the orbit $S^1\cdot c^m$ is an isolated
critical orbit of $E$.}

Since our aim is to prove the existence of more than one prime closed
geodesic for every Finsler metric on $S^2$, the condition (Iso) does not
restrict generality.

If $c$ has multiplicity $m$, then the subgroup
$$ \Z_m=\left\{\frac{k}{m}\,\,\mid\,\, 0\le k<m\right\} $$
of $S^1$ acts on $C_*(E,c)$. Quite generally, if $\Z_m$ acts on a set $H$,
we denote by $H^{\Z_m}$ the set of elements of $H$ fixed by $\Z_m$.

Our first aim is to prove

{\bf Proposition 3.2.} {\it If $c$ is a closed geodesic of multiplicity $m$
satisfying (Iso), then we have the following natural isomorphisms for all $q\in\Z$:
\bea
C_q(E,S^1\cdot c) &=& C_{q-1}(E,c)^{\Z_m}\oplus C_q(E,c)^{\Z_m},     \lb{2.10}\\
\overline{C}_q(E,c) &=& C_q(E,c)^{\Z_m}.  \lb{2.11}\eea}

The following lemmas will help prove Proposition 3.2.

{\bf Lemma 3.3.} {\it Let $\Lambda(k,a)\subseteq\Lambda$ be a
finite-dimensional approximation containing a closed geodesic $c$.
Let $D\subseteq\Lambda(k,a)$ be a hypersurface transverse to $S^1\cdot c$
at $c\in D$, and set $D^- = D\cap\Lambda(c)$. Then the inclusion
$D^-\cup\{c\}\rightarrow\Lambda(c)\cup\{c\}$ induces an isomorphism}
$$ H_*(D^-\cup\{c\}, D^-)= C_*(E,c). $$

{\bf Proof.} Using the deformation retraction $r_1=r(\cdot,1)$ defined
in (\ref{2.9}) we see that the inclusion induces an isomorphism
\be H_*\left((\Lambda(k,a)\cap\Lambda(c))\cup\{c\},
        \Lambda(k,a)\cap\Lambda(c)\right) = C_*(E,c). \lb{2.12}\ee
Next, we intend to deform a neighborhood $V\subseteq\Lm(k,a)$ of $c$
into $D$ without increasing energy. The energy non-increasing smooth map
$G:\Lambda(k,a)\times S^1\rightarrow\Lambda(k,a)$ defined by
$$   G(\gamma,s)=r_1(s\cdot\gamma) $$
is a submersion in a neighborhood $U$ of $(c,0)$ in $\Lambda(k,a)\times S^1$.

Since $\frac{\partial G}{\partial s}(c,0)$ is tangent to $S^1\cdot c$,
while $D$ is transverse to $S^1\cdot c$, we can find an open neighborhood
$V$ of $c$ in $\Lambda(k,a)$ and $\epsilon >0$ such that a smooth function
$\sigma:V\rightarrow (-\epsilon, \epsilon)$ is uniquely defined by
$$ G(\gamma,\sigma(\gamma))\in D.  $$
Now we define $h:V\times[0,1]\rightarrow\Lm(k,a)$ by
$$ h(\ga,t) = G(\ga,t\sg(\ga)) = r_1((t\sg(\ga))\cdot\ga).  $$
Then we have $h_0 =\mbox{id}_V$, $h_1(V)\subseteq D$, $h(\ga,t)=\ga$ for
every $\ga\in D\cap V$, $t\in[0,1]$, and $E\circ h\le E$. Using the homotopy
$h$ and excision one can see that the inclusion
$D^-\cup\{c\}\rightarrow(\Lambda(k,a)\cap\Lambda(c))\cup\{c\}$ induces an
isomorphism
\be H_*(D^- \cup \{c\}, D^-) =
      H_*((\Lm(k,a)\cap\Lm(c))\cup\{c\},\Lm(k,a)\cap\Lm(c)). \lb{2.13}
\ee
Now (\ref{2.12}) and (\ref{2.13}) imply our claim. \hb

We need the following variants of Lemma 3.3 that involve the isotropy group of $c$.

{\bf Lemma 3.4.} {\it Let $c$ be a closed geodesic of multiplicity $m\ge 2$ and
$\Lambda(j,a)\subseteq\Lambda$ a finite-dimensional approximation containing $c$
and such that $m$ divides $j$. Let $D\subseteq\Lambda(j,a)$ be a $\Z_m$-invariant
hypersurface transverse to $S^1\cdot c$ at $c\in D$, and set $D^-=D\cap\Lambda(c)$.
Then the inclusion $D^-\cup\{c\}/\Z_m\rightarrow\Lambda(c)\cup\{c\}/\Z_m$ induces
an isomorphism
$$ H_*(D^-\cup\{c\}/\Z_m, D^-/\Z_m)
\rightarrow H_*(\Lambda(c)\cup\{c\}/\Z_m, \Lambda(c)/\Z_m).  $$
Moreover, let $\Gamma\simeq\Z_m$ act on $(\Lambda(c)\cup\{c\})\times S^1$ by
\be \frac{k}{m}(\ga,s)=\left(\frac{k}{m}\cdot\ga, s-\frac{k}{m}\right),
      \qquad  {\it for}\;\;k=0,\ldots,m-1.  \lb{2.14}\ee
Then the inclusion
$((D^-\cup\{c\})\times S^1)/\Gamma \rightarrow (\Lambda(c)\times S^1)/\Gamma$
induces an isomorphism
$$ H_* (((D^-\cup\{c\})\times S^1)/\Gamma, (D^-\times S^1)/\Gamma)\rightarrow
H_* (((\Lambda(c)\cup\{c\})\times S^1)/\Gamma, (\Lambda(c)\times S^1)/\Gamma).
$$}

{\bf Remark.} If one applies the exponential map of $\Lambda(j,a)$ to a
small neighborhood of the origin in the normal space at $c$ of $S^1\cdot c$
in $\Lambda(j,a)$, then one obtains a hypersurface $D$ satisfying the
assumptions made in Lemma 3.4.

{\bf Proof.} All the constructions in the proof of Lemma 3.3 are
$\Z_m$-equivariant. In particular, one can choose the neighborhood $V$ of
$c$ in $\Lambda(j,a)$ to be $\Z_m$-invariant, and then the homotopy $h$
satisfies
$h\left(\frac{k}{m}\cdot\gamma, t\right) = \frac{k}{m}\cdot h(\gamma,t)$ for
all $k=0,\ldots,m-1$, $\gamma\in V$ and $t\in[0,1]$. This implies the first
statement of Lemma 3.4. The proof for the second statement is analogous. Here
one defines a homotopy
$\td{h}:((V\times S^1)/\Ga)\times[0,1]\rightarrow(\Lm(k,a)\times S^1)/\Ga$ by
$$ \td{h}([\ga,s],t)=[h(\ga,t),s], $$
where the square brackets denote the $\Ga$-orbits. \hb

{\bf Lemma 3.5.} {\it Assume that the closed geodesic $c\in\Lm(k,a)$ satisfies
(Iso) and $D\subseteq\Lm(k,a)$ is a hypersurface transverse to $S^1\cdot c$
at $c\in D$. Then $c$ is an isolated critical point of $E|_D$. }

{\bf Proof.} It suffices to prove that for $\ga\in D\setminus\{c\}$ close to
$c$, the hyperplane $T_{\ga}D$ is not contained in the kernel of $E'(\ga)$.
Consider the curve $\Gamma_{\ga}:S^1\to \Lm(k,a)$ defined by
$$  \Gamma_{\ga}(s) = r_1(s\cdot\ga), $$
where $r_1$ denotes the retraction $\Lm^{a_-}\to \Lm(k,a)$ used in the
proof of Lemma 3.3. Since we have
$E\circ\Gamma_{\ga}(s)\le E\circ\Gamma_{\ga}(1)$ for all $s\in S^1$,
we see that
\be \Gamma'_{\ga}(1)\in\ker E'(\Gamma_{\ga}(1))=\ker E'(\ga). \lb{2.15}\ee
Since $\Gamma'_c(1)$ is tangent to $S^1\cdot c$ and $D$ is transverse to
$S^1\cdot c$ in $c$, we have $\Gamma'_c(1)\notin T_c D$. By continuity we
see that
\be  \Gamma'_{\ga}(1)\not\in T_{\ga}D. \lb{2.16}\ee
if $\gamma\in D$ is close to $c$. Since $S^1\cdot c$ is an isolated critical
orbit, all $\ga\in D\setminus\{c\}$ close to $c$ are regular points of
$E|_{\Lm(k,a)}$. Therefore both $\ker E'(\ga)\cap T_{\ga}\Lm(k,a)$ and
$T_{\ga}D$ are codimension $1$ subspaces in $T_{\ga}\Lm(k,a)$. By (\ref{2.15})
and (\ref{2.16}), they are different. Thus there exists $\xi\in T_{\ga}D$
such that $E'(\ga)(\xi)\not= 0$. \hb

{\bf Lemma 3.6.} {\it Let $D$ be a finite-dimensional Riemannian manifold on
which $\Z_m\subseteq S^1$ acts by isometries. Let $E:D\rightarrow\R$ be smooth
and $\Z_m$-invariant, and suppose $c\in D$ is the only critical point of $E$.
Set $D^- = E^{-1}((-\infty, E(c)))$. Then the transfer homomorphism
$$ H_*(D^-\cup\{c\}/\Z_m, D^-/\Z_m) \rightarrow H_*(D^-\cup\{c\}, D^-)^{\Z_m} $$
is an isomorphism. }

{\bf Proof.} This depends on the fact that we use homology with rational
coefficients. If we would use \u{C}ech homology instead of singular homology,
Lemma 3.6 would directly follow from \cite{Bre1}, Theorem III.7.2. Since
\u{C}ech and singular homology coincide for triangulable pairs, it is sufficient
to show that, after appropriate excision, $(D^-\cup\{c\}/\Z_m, D^-/\Z_m)$ is
homotopy equivalent to a triangulable pair. Note that, by \cite{Cha1}, Theorem
I.7.8, we can perturb $E$ in an arbitrarily small neighborhood of $c$ to a
$\Z_m$-invariant function with only non-degenerate critical points. Then
$\Z_m$-equivariant Morse theory, cf.~\cite{Was1}, can be applied to finish
the proof. A different argument proving Lemma 3.6 is contained in \cite{Rad2},
Section 6.3. Here, H.-B. Rademacher uses relative $\Z_m$-CW-complexes. \hb

{\bf Proof of Proposition 3.2.} Choose a tubular neighborhood $W$ of the
circle $S^1\cdot c$ in $\Lambda$ with fibers $W_{s\cdot c}=s\cdot W_c$ over
$s\cdot c\in S^1\cdot c$, cf. \cite{Kli2}, Lemma 2.2.8. As on pp.502-503
of \cite{GrM2} we see that the map
$$ W_c\times S^1\rightarrow W, \,\, (\gamma, s)\mapsto s\cdot\gamma, $$
is a normal covering with group of covering transformations
$\Gamma\simeq\Z_m$ operating by (\ref{2.14}). This together with excision
provides an isomorphism
\be H_*(((W^-_c\cup\{c\})\times S^1)/\Gamma,
     (W^-_c\times S^1)/\Gamma) = C_*(E, S^1\cdot c).  \lb{2.17}\ee
Here and below we set $W^-=W\cap \Lm(c)$ and $W_c^-=W_c\cap\Lm(c)$.
Now choose a finite-dimensional approximation $\Lm(j,a)\subseteq\Lm$
and a hypersurface $D\subseteq W_c$ in $\Lambda(j,a)$ as in Lemma 3.4.
Since some neighborhood of $c$ in $\Lambda$ is $E$-equivariantly
homeomorphic to $W_c\times(-\epsilon,\epsilon)$ with
$\epsilon\in\left(0, \frac{1}{2m}\right)$, we can use excision to obtain
an isomorphism
\be H_*\left(\left((W^-_c\cup\{c\})\times S^1\right)/\Gamma,
(W^-_c\times S^1)/\Gamma\right)
  = H_*\left(\left((\Lambda(c)\cup\{c\})\times S^1\right)/\Gamma,
  (\Lambda(c)\times S^1)/\Gamma\right)   \lb{2.18}
\ee
Then (\ref{2.17}), (\ref{2.18}) and Lemma 3.4 provide an isomorphism
\be H_*\left(\left((D^-\cup\{c\})\times S^1\right)/\Gamma,
    (D^-\times S^1)/\Gamma \right) = C_*(E,S^1\cdot c). \lb{2.19}\ee
Now the proof of Lemma 2.6 shows that the transfer is an isomorphism
\be
H_*\left((D^-\cup\{c\})\times S^1,D^-\times S^1\right)^{\Gamma}
 = H_*\left(\left((D^-\cup\{c\})\times S^1\right)/\Gamma,
    (D^-\times S^1)/\Gamma\right)  \lb{2.20}\ee
Using the K\"unneth formula we see that (\ref{2.19}) and (\ref{2.20})
provide an isomorphism
$$ C_q(E, S^1\cdot c) =
  H_{q-1}(D^-\cup\{c\}, D^-)^{\Z_m}\oplus H_q(D^-\cup\{c\}, D^-)^{\Z_m} $$
for all $q\in\Z$. By Lemma 3.3 this implies (\ref{2.10}).

Finally, the pair $(W^-\cup S^1\cdot c/S^1, W^-/S^1)$ is obviously
homeomorphic to the pair $(W^-_c\cup\{c\}/\Z_m, W^-_c/\Z_m)$.
Using Lemma 2.4 and similar arguments as before, we obtain an isomorphism
$$  H_*(D^-\cup\{c\}/\Z_m, D^-/\Z_m)=\overline{C}_*(E,c)  $$
Now we can apply Lemma 3.3 and Lemma 3.6 to obtain (\ref{2.11}). \hb

We will now apply the results by D. Gromoll and W. Meyer \cite{GrM1} to a
given closed geodesic $c$ satisfying (Iso). If $m=m(c)$ is the multiplicity
of $c$, we choose a finite-dimensinal approximation $\Lm(k,a)\subseteq\Lm$
containing $c$ such that $m$ divides $k$. Then the isotropy subgroup
$\Z_m\subseteq S^1$ of $c$ acts on $\Lm(k,a)$ by isometries.
Let $D$ be a $\Z_m$-invariant local hypersurface transverse to $S^1\cdot c$
at $c\in D$. Such $D$ can be obtained by applying the exponential map of
$\Lambda(k,a)$ at $c$ to the normal space to $S^1\cdot c$ at $c$. We let
\be T_c D=V_+\oplus V_-\oplus V_0    \lb{2.21}\ee
denote the orthogonal decomposition of $T_c D$ into the positive, negative and
null eigenspace of the endomorphism of $T_c D$ associated to $(E|_D)''(c)$ by
the Riemannian metric. In particular, we have $\dim V_- = i(c)$ and
$\dim V_0=\nu(c)$. According to \cite{GrM1}, Lemma 1, there exist balls
$B_+\subseteq V_+$, $B_-\subseteq V_-$, and $B_0\subseteq V_0$ centered at the
origins and a diffeomorphism
\be \psi:B_+\times B_-\times B_0\rightarrow
   \psi(B_+\times B_-\times B_0)\subseteq D    \lb{2.22}\ee
such that $\psi(0,0,0)=c$, $\psi_{*(0,0,0)}$ preserves the splitting
(\ref{2.21}), and
\be  E\circ\psi(x_+, x_-, x_0) = |x_+|^2 - |x_-|^2 + f(x_0), \lb{2.23}\ee
where $f:B_0\rightarrow\R$ satisfies $f'(0)=0$ and $f''(0)=0$.
Since the $\Z_m$-action is isometric and preserves $E$, the differential
$\left(\frac{1}{m}|_D\right)_{*c}$ preserves the splitting (\ref{2.21}).
It follows from the construction of $\psi$ that $\psi$ is equivariant with
respect to the $\Z_m$-action, i.e.,
\be \frac{1}{m}\circ\psi=\psi\circ\left(\frac{1}{m}|_D\right)_{*c},
     \lb{2.24}\ee
cf. \cite{GrM2}, p.501. As usual, we call
$$ U=\{\psi(0,x_-, 0)\;|\;x_-\in B_-\} $$
a local negative disk at $c$, and
$$ N=\{\psi(0,0,x_0)\;|\;x_0 \in B_0\} $$
a local characteristic manifold at $c$. By (\ref{2.24}), local negative
disks and local characteristic manifolds are $\Z_m$-invariant.

It follows from Lemma 3.5 and (\ref{2.23}) that $c$ is an isolated critical
point of $E|_N$. We set $N^-=N\cap\Lm(c)$, $U^-=U\cap\Lm(c)=U\setminus\{c\}$
and $D^-=D\cap\Lm(c)$. Using (\ref{2.23}) and the fact that $c$ is an
isolated critical point of $E|_N$ and the K\"unneth formula one concludes
\be H_*(D^-\cup\{c\},D^-)=
      H_*(U^-\cup\{c\},U^-) \otimes H_*(N^-\cup\{c\},N^-), \lb{2.25}\ee
where
\be
H_q(U^-\cup\{c\}, U^-) = H_q(U, U\setminus\{c\})
    = \left\{\matrix{\Q, & {\rm if\;}q=i(c), \cr
                      0, & {\rm otherwise}, \cr}\right.  \lb{2.26}\ee
cf. \cite{Rad2}, Lemma 6.4 and its proof, or \cite{GrM1}, Lemma 6.

Using Proposition 3.2, Lemma 3.3, (\ref{2.25}) and (\ref{2.26}), we
obtain the following version of the Gromoll-Meyer Shifting Lemma:

{\bf Proposition 3.7.} {\it Suppose $c$ is a closed geodesic of multiplicity
$m(c)=m$ satisfying (Iso). Let $U$ be a local negative disk at $c$ and let
$N$ be a local characteristic manifold at $c$. Then for $q\in \Z$ there holds }
\bea
C_q(E,S^1\cdot c)
&=& \left(H_{i(c)}(U^-\cup\{c\},U^-)
      \otimes H_{q-i(c)}(N^-\cup\{c\},N^-)\right)^{\Z_m} \nn\\
&&\qquad\quad \oplus \left(H_{i(c)}(U^-\cup\{c\},U^-)
  \otimes H_{q-1-i(c)}(N^-\cup\{c\},N^-)\right)^{\Z_m}. \nn
\eea

In order to obtain a more explicit version of the formula in Proposition 3.7
one needs to know if a generator of the $\Z_m$-action on $U$ reverses
orientation or not. This has been investigated by A.S. Svarc \cite{Sva1}, see
also \cite{Kli1}, \cite{Kli2}, Lemma 4.1.4, and \cite{Rad2}, Section 6.3.

We introduce the following notation. If the group $\Z_m$ acts linearly on
a vector space $H$ and if $T$ is a generator of $\Z_m$, we let
$H^{\Z_m,1}=H^{\Z_m}$ denote the eigenspace of $T$ corresponding to $1$,
while $H^{\Z_m,-1}$ denotes the eigenspace of $T$ corresponding to $-1$.
This is independent of the choice of the generator $T$ of $\Z_m$. If $m$ is
odd then $H^{\Z_m,-1}=\{0\}$.

{\bf Proposition 3.8.} {\it Let $c$ be a prime closed geodesic satisfying
(Iso) and let $m\in\N$.
\begin{itemize}
\item[(i)]
If $\nu(c^m)=0$, then
$$ C_q(E,S^1\cdot c^m)=
  \left\{\matrix{\Q, & {\rm if\;}i(c^m) - i(c)\in 2\Z,
                         {\it and\;}q \in\{i(c^m), i(c^m)+1\} \cr
                 0,  & {\it otherwise}. \cr}\right.  $$
\item[(ii)]
If $\nu(c^m)>0$, let $N_{c^m}$ be a local characteristic manifold at
$c^m$ and $N^-_{c^m}=N_{c^m}\cap\Lambda(c^m)$. We set
$\ep(c^m)=(-1)^{i(c^m)-i(c)}$. Then we have
\newline
$C_q(E,S^1\cdot c^m)=
  H_{q-i(c^m)}(N^-_{c^m}\cup\{c^m\},N^-_{c^m})^{\Z_m,\ep(c^m)}\oplus
  H_{q-1-i(c^m)}(N^-_{c^m}\cup\{c^m\},N^-_{c^m})^{\Z_m,\ep(c^m)}$.
\end{itemize}}

{\bf Proof.} Let $U$ denote a local negative disk at $c^m$. The question of
whether a generator $T$ of the $\Z_m$-action on $U$ acts as multiplication by $+1$
or $-1$ on $H_{i(c^m)}(U,U\setminus\{c^m\})=\Q$ reduces to the question whether
$\det(D(T|_U)_{c^m})$ is positive or negative. Let $a_+$ and $a_-$ denote
the dimensions of the $(+1)$-eigenspace and the $(-1)$-eigenspace of $D(T|_U)_{c^m}$
respectively. Since $T$ is an isometry we see that
$$  a_+ + a_- = \dim(U)=i(c^m) \quad\mod\; 2, $$
and
$$ \det(D(T|_U)_{c^m}) = (-1)^{a_-}.  $$
Since $DT_{c^m}(\xi)=\xi$ for some $\xi\in T_{c^m}\Lm$ if and only if
$\xi=D\phi_m(\tilde{\xi})$ for some $\tilde{\xi}\in T_c\Lm$, we conclude
from (\ref{2.6}) that $a_+ = i(c)$ and hence
$$  a_- = i(c^m)-i(c)\quad\mod \;\;2.  $$
This implies
$$  \det(D(T|_U)_{c^m})=\epsilon(c^m).  $$
Now our claim follows from Proposition 3.7. \hb

{\bf Definition 3.9.} {\it  \begin{itemize}
\item[(i)] Suppose $c$ is a closed geodesic such that $\nu(c^n)>0$ for
some $n\in\N$. Then we set
$$ n_c=n(c)=\min\{n\in\N\,|\,\nu(c^n)>0\}. $$
\item[(ii)] Suppose $c$ is a closed geodesic of multiplicity $m(c)=m$
satisfying (Iso). If $N$ is a local characteristic manifold at $c$,
$N^- = N\cap\Lm(c)$ and $j\in\N\cup\{0\}$, we define
\bea
k_j(c) &=& \dim H_j(N^-\cup\{c\},N^-), \nn\\
\hat{k}_j(c) &=& \dim H_j(N^-\cup\{c\},N^-)^{\Z_m}. \nn
\eea
\end{itemize}}

Note that Lemma 3.3 and (\ref{2.25}) imply that the numbers $k_j(c)$
and $\hat{k}_j(c)$ are independent of the choice of $N$. Moreover, we
obviously have
$$  0\le\hat{k}_j(c)\le k_j(c).  $$
Note that the finiteness of $k_j(c)$ follows from Lemma 2 of \cite{GrM1}.
Since $\Z_m$ fixes $c$, we obviously have
\be k_0(c)=\hat{k}_0(c).  \lb{2.27}\ee
Finally, if $c$ is non-degenerate, i.e., if $\nu(c)=0$, then
$k_0(c)=\hat{k}_0(c)=1$, while $k_j(c)=0$ for $j>0$.

The following facts will be useful.

{\bf Lemma 3.10.} {\it Let $c$ be a closed geodesic satisfying (Iso),
and let $N$ be a local characteristic manifold at $c$. Then the following
is true.
\begin{itemize}
\item[(i)] The closed geodesic $c$ is a strict local minimum of $E|_N$ if
and only if $k_0(c)\not= 0$. In particular, $k_0(c)\not= 0$ implies
$k_0(c)=1$ and $k_j(c)=0$ for $j>0$.
\item[(ii)] The closed geodesic $c$ is a strict local maximum of $E|_N$ if
and only if $k_{\nu(c)}(c)\not= 0$. In particular, $k_{\nu(c)}(c)\not= 0$
implies $k_{\nu(c)}(c)=1$ and $k_j(c)=0$ for $j\not=\nu(c)$.
\end{itemize}}

{\bf Remark.} Note that $c$ is an isolated critical point of $E|_N$.
Hence, if $c$ is a local minimum of $E|_N$, then this is strictly so, and
similarly for the condition "local maximum".

{\bf Proof.} (i) Assume $k_0(c)\not= 0$, i.e., $H_0(N^-\cup\{c\},N^-)\not= 0$.
Since a connected open neighborhood $V$ of $c$ in $N$ can be continuously
deformed into $N^-\cup\{c\}$ in an energy non-increasing manner, cf.
\cite{Cha1}, Theorem I.3.2., we conclude that $V^- = V\cap\Lm(c)=\emptyset$.
Hence $c$ is a strict local minimum of $E|_N$.

(ii) This is proved in \cite{Hin1}, p.256. \hb

We now come to Gromoll-Meyer's crucial result on the type numbers $k_j(c^m)$
of an iterated prime closed geodesic $c$, cf. \cite{GrM2}, Theorem 3. Such
a study for Lagrangian systems was carried out in \cite{Lon2} (cf. also
Section 14.3 of \cite{Lon3}). We obtain a similar result for the dimensions
$\hat{k}_j(c^m)$ of the $\Z_m$-invariant part too.

{\bf Theorem 3.11.} {\it Let $c$ be a prime closed geodesic in a Finsler
manifold. Suppose $c$ satisfies (Iso), $m$, $n$, $p$ are integers and
$m=np$. If the nullities of $c^m$ and $c^n$ satisfy
$$  \nu(c^m) = \nu(c^n),  $$
then
$$ k_j(c^m) = k_j(c^n)  \quad {\rm and}\quad
      \hat{k}_j(c^m) = \hat{k}_j(c^n) $$
for all $j\in\N\cup\{0\}$. }

{\bf Proof.} We choose finite-dimensional approximations $\Lm(k,a)$ containing
$c^n$ and $\Lm(kp,p^2a)$ containing $c^m$ and a characteristic manifold
$N\subseteq\Lm(k,a)$ at $c^n$. Note that the iteration map $\phi_p$ defined
by (\ref{2.3}) maps $\Lm(k,a)$ diffeomorphically to a submanifold of
$\Lm(kp,p^2a)$. Hence $\phi_p(N)$ is a submanifold of $\Lm(kp,p^2a)$ transverse
to $S^1\cdot c^m$ whose tangent space at $c^m$ is contained in the null space
of $E''(c^m)$, cf. Lemma 3.1. Note moreover that
$\dim \phi_p(N)=\dim N=\nu(c^n)$ and $\nu(c^n)=\nu(c^m)$ by assumption. Arguing
as in the proof of Theorem 3 in \cite{GrM2}, we can now invoke Lemma 7 from
\cite{GrM1} to conclude that $\phi_p(N)$ is a characteristic manifold at $c^m$.
Since $E\circ\phi_p=p^2 E$, this implies that $k_j(c^m)=k_j(c^n)$ for all $j$.
If $T_{\frac{1}{n}}$ and $T_{\frac{1}{m}}$ denote the actions of $\frac{1}{n}$
and $\frac{1}{m}$ on $N$ and $\phi_p(N)$, respectively, then
$$  T_{\frac{1}{m}}\circ\phi_p = \phi_p\circ T_{\frac{1}{n}}, $$
cf. \cite{Rad2}, p.67. Hence
$$ (\phi_p)_* : H_*(N^-\cup\{c^n\},N^-)
       \rightarrow H_*(\phi_p(N)^-\cup\{c^m\},\phi_p(N)^-)  $$
is an isomorphism conjugating generators of the $\Z_n$-action on
$H_*(N^-\cup\{c^n\},N^-)$ and the
$\Z_m$-action on $H_*(\phi_p(N)^-\cup\{c^m\},\phi_p(N)^-)$.
This implies $\hat{k}_j(c^n)=\hat{k}_j(c^m)$ for all $j$. \hb

The following result on the critical modules $\overline{C}_*(E,c)$ in the
quotient $\overline{\Lambda}=\Lambda/S^1$, cf. \cite{Rad2}, Satz 6.11, will
be used in Section 5.

{\bf Proposition 3.12.} {\it Let $c$ be a prime closed geodesic satisfying
$(Iso)$ and let $m\in\N$, $q\in\N\cup\{0\}$. Let $N$ be a characteristic
manifold at $c^m$, $N^- =N\cap\Lm(c^m)$. Then we have }
$$ \ol{C}_q(E,c^m) = H_{q-i(c^m)}(N^-\cup\{c^m\},N^-)^{\Z_{m},\ep(c^m)}. $$

{\bf Proof.} Using (\ref{2.11}), Lemma 3.3, and (\ref{2.25}) we conclude
that
\be \overline{C}_q(E, c^m)=
\left(H_{i(c^m)}(U^-\cup\{c^m\},U^-)\otimes H_{q-i(c^m)}(N^-\cup\{c^m\},N^-))
    \right)^{\Z_m}, \lb{2.28}
\ee
if $U$ is a local negative disk at $c^m$. Since generators of $\Z_m$ act
on $H_{i(c^m)}(U^-\cup\{c^m\},U^-)=\Q$ through multiplication by
$\ep(c^m)=(-1)^{i(c^m)-i(c)}$, cf. the proof of Proposition 3.8, our
claim follows from (\ref{2.28}). \hb

In order to relate the critical $\Q$-modules $C_*(E, S^1\cdot c)$ of closed
geodesics $c$ to the homology of the loop space $\Lambda$, we will use the
following fact.

{\bf Proposition 3.13.} {\it Suppose $u\in(a,b)$ is the only critical value
of $E$ in the interval $[a,b]$, and the critical set
$$ C=\{c\,|\, c \mbox{ is a closed geodesic with $E(c)=u$}\} $$
is the disjoint union of finitely many critical orbits
$$ C=\bigcup\limits^q_{i=1}S^1\cdot c_i. $$
Then there is an isomorphism
\be \bigoplus\limits^q_{i=1}C_*(E, S^1\cdot c_i)=H_*(\Lambda^b, \Lambda^a).
  \lb{2.29}\ee}

{\bf Proof.} We choose a finite-dimensional approximation
$\tilde{\Lambda}=\Lambda(\tilde{k}, \tilde{a})$ with $\tilde{a}>b$, and we
set $\tilde{\Lambda}^b:=\Lambda^b\cap\tilde{\Lambda}$,
$\tilde{E}=E|_{\tilde{\Lambda}}$ etc. Using the energy non-increasing
deformation retraction from $\Lambda^{\tilde{a}_-}$ onto $\tilde{\Lambda}$,
one sees that it suffices to prove
\be \bigoplus\limits^q_{i=1}C_*(\tilde{E}, S^1\cdot c_i)
    =H_*(\tilde{\Lambda}^b, \tilde{\Lambda}^a),   \lb{2.30}\ee
where $C_*(\td{E},S^1\cdot c_i)=
H_*(\td{\Lm}^{u_-}\cup S^1\cdot c_i, \td{\Lm}^{u_-})$.

Note that
\be H_*(\tilde{\Lambda}^{u_-}\cup C, \tilde{\Lambda}^{u_-})
    = \bigoplus\limits^q_{i=1} C_*(\tilde{E}, S^1\cdot c_i).  \lb{2.31}\ee

Since there are no critical values of $\tilde{E}$ in the interval $[a,u)$,
the flow of -grad $\tilde{E}$ induces a strong deformation retraction of
$\tilde{\Lambda}^{u_-}$ onto $\tilde{\Lambda}^a$. This implies
$$ H_*(\tilde{\Lambda}^b, \tilde{\Lambda}^a)
    = H_*(\tilde{\Lambda}^b, \tilde{\Lambda}^{u_-}). $$
Next we choose disjoint tubular neighborhoods $W_i\subseteq\tilde{\Lambda}^b$
of the critical orbits $S^1\cdot c_i$. By Lemma 3.5 we can assume that the
orthogonal projection -grad $\tilde{E}^{\top}$ of -grad $\tilde{E}$ to the
tangent spaces of the fibers of the tubular neighborhoods vanishes only on
the critical set $C$.

We choose a smooth function $\lambda:\tilde{\Lambda}^b\rightarrow[0,1]$ with
support in $\bigcup\limits^q_{i=1} W_i$ and such that $\lm=1$ holds in a
neighborhood of $C$. Now we consider the vector field
$$ X=(1-\lm)(-\mbox{grad }\td{E})+\lm(-\mbox{grad }\td{E}^{\top}) $$
on $\tilde{\Lambda}^b$. Since the restrictions of $\tilde{E}$ to the fibers
of the tubular neighborhoods have only one critical point, one easily sees
that a flow line of $X$ starting in $\td{\Lm}^b\bs\tilde{\Lm}^{u_-}$ either
reaches $\td{\Lm}^{u_-}$ in finite time or converges to some single point in
$C\subset \td{\Lm}^b$. This allows us to define a homotopy
$$ H:\td{\Lm}^b\times[0,1]\rightarrow\td{\Lm}^b  $$
such that
$$ H(\tilde{\Lambda}^{u_-}\times[0,1])\subseteq\tilde{\Lambda}^{u_-}, $$
and $H_1=H(\cdot,1)$ satisfies
$$ H_1(\tilde{\Lambda}^b)\subseteq\tilde{\Lambda}^{u_-}\cup C,  $$
see e.g. \cite{Cha1}, Theorem I.3.2. This implies
$$ H_*(\td{\Lm}^b, \td{\Lm}^{u_-}) = H_*(\td{\Lm}^{u_-}\cup C, \td{\Lm}^{u_-}).
$$
Using this and (\ref{2.31}) we conclude that (\ref{2.30}) is true. \hb

Suppose $(M,F)$ is a compact Finsler manifold that has only $q$ prime closed
geodesics $c_j$ for $1\le j\le q$. Then the Morse type numbers $M_k$ for
$k\in\N\cup\{0\}$ are defined by
$$ M_k = \sum_{1\le j\le q \atop m\ge 1}\dim C_k(E, S^1\cdot c^m_j). $$
Using Proposition 3.13 we can prove the Morse inequalities in the standard
fashion, see e.g. \cite{Kli2}, Theorem 2.4.12, or \cite{Cha1}, Theorem I.4.3.
Let $b_k=b_k(\Lambda, \Lambda^0)=\dim H_k(\Lambda, \Lambda^0)$ denote the
relative Betti numbers of the pair $(\Lm,\Lm^0)$ with coefficients in $\Q$.

{\bf Theorem 3.14.} {\it Let $(M,F)$ be a compact Finsler manifold with only
finitely many prime closed geodesics. Then for every integer $k\ge 0$ there
holds
\be  M_k\ge b_k=b_k(\Lambda, \Lambda^0)  \lb{2.32}\ee
and
\be\matrix{
& M_k - M_{k-1} + M_{k-2} -  \cdots + (-1)^{k-1} M_1 + (-1)^k M_0   \cr
& \qquad \ge b_k - b_{k-1} + b_{k-2}  - \cdots (-1)^{k-1} b_k + (-1)^k b_0. \cr}
\lb{2.33}\ee}

\setcounter{equation}{0}
\section{Classification of closed geodesics on $S^2$}

Let $c$ be a closed geodesic on a Finsler sphere $S^2=(S^2,F)$. Denote the
linearized Poincar\'e map of $c$ by $P_c:\R^2\to\R^2$. Because the index
iteration formulae in \cite{Lon3} work for Morse index of iterated closed
geodesics, the iteration formula of Morse indices of $c$ must be one of
the following nine cases by Theorems 8.1.4 to 8.1.7 of \cite{Lon3}. Here we
use the notation from Section 8.1 of \cite{Lon3}.

{\bf Case CG-1.} {\it $P_c$ is conjugate to a matrix
$\left(\matrix{ 1 & b\cr
                0 & 1\cr}\right)$ for some $b>0$.}

In this case, by $1^{\circ}$ of Theorem 8.1.4 of \cite{Lon3}, we have
$i(c)=2p-1$ for some $p\in\N$, and
\be i(c^m)=2mp-1, \quad \nu(c^m)=1, \qquad {\rm for\;all}\; m\ge 1. \lb{3.1}\ee

{\bf Case CG-2.} {\it $P_c=I_2$, the $2\times 2$ identity matrix}.

In this case by $2^{\circ}$ of Theorem 8.1.4 of \cite{Lon3}, we have
$i(c)=2p-1$ for some $p\in\N$, and
\be i(c^m)=2mp-1, \quad \nu(c^m)=2, \qquad {\rm for\;all}\; m\ge 1. \lb{3.2}\ee

{\bf Case CG-3.} {\it $P_c$ is conjugate to a matrix
$\left(\matrix{ 1 & -b\cr
                0 & 1\cr}\right)$ for some $b>0$.}

In this case by $3^{\circ}$ of Theorem 8.1.4 of \cite{Lon3}, we have
$i(c)=2p$ for some $p\in\N\cup\{0\}$, and
\be i(c^m)=2mp, \quad \nu(c^m)=1, \qquad {\rm for\;all}\; m\ge 1. \lb{3.3}\ee

{\bf Case CG-4.} {\it $P_c$ is conjugate to a matrix
$\left(\matrix{ -1 & -b\cr
                 0 & -1\cr}\right)$ for some $b>0$}.

In this case by $1^{\circ}$ of Theorem 8.1.5 of \cite{Lon3}, we have
$i(c)=2p+1$ for some $p\in\N\cup\{0\}$, and
\be i(c^m)=m(2p+1)-\frac{1+(-1)^m}{2}, \quad \nu(c^m)= \frac{1+(-1)^m}{2},
          \qquad {\rm for\;all}\; m\ge 1. \lb{3.4}\ee

{\bf Case CG-5.} {\it $P_c=-I_2$}.

In this case by $2^{\circ}$ of Theorem 8.1.5 of \cite{Lon3}, we have
$i(c)=2p+1$ for some $p\in\N\cup\{0\}$, and
\be i(c^m)=m(2p+1)-\frac{1+(-1)^m}{2}, \quad \nu(c^m)= 1+(-1)^m,
          \qquad {\rm for\;all}\; m\ge 1. \lb{3.5}\ee

{\bf Case CG-6.} {\it $P_c$ is conjugate to a matrix
$\left(\matrix{ -1 & b\cr
                 0 & -1\cr}\right)$ for some $b>0$.}

In this case by $3^{\circ}$ of Theorem 8.1.5 of \cite{Lon3}, we have
$i(c)=2p+1$ for some $p\in\N\cup\{0\}$, and
\be i(c^m)=m(2p+1), \quad \nu(c^m)= \frac{1+(-1)^m}{2},
          \qquad {\rm for\;all}\; m\ge 1. \lb{3.6}\ee

{\bf Case CG-7.} {\it $P_c$ is rationally elliptic, i.e., $P_c$ is
conjugate to some rotation matrix
$R(\th)=\left(\matrix{\cos\th & -\sin\th\cr
                      \sin\th & \cos\th\cr}\right)$ with some
$\th\in (0,\pi)\cup (\pi,2\pi)$ and $\th/\pi\in\Q$}.

In this case by Theorem 8.1.7 of \cite{Lon3}, we have $i(c)=2p+1$
for some $p\in \N\cup\{0\}$, and
\be
i(c^m) = \left\{\matrix{
   2mp + 2[\frac{m\th}{2\pi}] + 1, &\quad \nu(c^m)=0,
          &\qquad {\rm if}\;\;m\th\not= 0 \;\;\mod\;\;2\pi, \cr
   2mp + 2[\frac{m\th}{2\pi}] - 1, &\quad \nu(c^m)=2,
          &\qquad {\rm if}\;\;m\th = 0 \;\;\mod\;\;2\pi. \cr}\right.
\lb{3.7}\ee

{\bf Case CG-8.} {\it $P_c$ is irrationally elliptic, i.e., $P_c$ is
conjugate to some rotation matrix $R(\th)$ with some
$\th\in (0,\pi)\cup (\pi,2\pi)$ and $\th/\pi\not\in\Q$}.

In this case by Theorem 8.1.7 of \cite{Lon3}, we have $i(c)=2p+1$ for
some $p\in \N\cup\{0\}$, and
\be i(c^m) = 2mp + 2[\frac{m\th}{2\pi}] + 1, \quad \nu(c^m)=0,
          \qquad {\rm for\;all}\; m\ge 1. \lb{3.8}\ee

{\bf Case CG-9.} {\it $P_c$ is hyperbolic, i.e., $P_c$ is conjugate
to the matrix $\left(\matrix{ b & 0\cr
                  0 & 1/b\cr}\right)$ for some $b>0$ or $b<0$.}

In this case, by Theorem 8.1.6 of \cite{Lon3}, we have $i(c)=p$
for some $p\in \N\cup\{0\}$, and
\be i(c^m) = mp, \quad \nu(c^m)=0, \qquad {\rm for\;all}\; m\ge 1. \lb{3.9}\ee

It is well known that if all iterations $c^m$ of a closed geodesic
$c$ are non-degenerate, $c$ must be hyperbolic or irrationally
elliptic, i.e., $P_c$ is of the class CG-8 or CG-9. In this case, $c$
is called {\it non-degenerate}.

{\bf Remark 4.1.} We should remind the readers that for a closed
geodesic $c:\R/(\tau\Z)\to M$ on a Finsler surface $M$, the linearized
Poincar\'e map $P_c\in \Sp(2)$ is given by
$\left(\matrix{x(\tau) & y(\tau) \cr
               \dot{x}(\tau) & \dot{y}(\tau) \cr}\right)$ in Section 3.4
of \cite{Kli1} and \cite{Kli2}. But in the notations here as well as in
\cite{Lon1} and \cite{Lon3}, the matrix $P_c$ is given by
$\left(\matrix{\dot{y}(\tau) & \dot{x}(\tau) \cr
               y(\tau) & x(\tau) \cr}\right)$.

\setcounter{equation}{0}
\section{Rationally elliptic degenerate saddle closed geodesics on $S^2$}

In this section, we study a particular type of closed geodesics of class CG-7.

{\bf Definition 5.1.} {\it Let $c$ be a prime
closed geodesic on a Finsler $2$-sphere that is rationally elliptic, i.e.,
its linearized Poincar\'e map $P_c$ is of class CG-7 with rotation angle
$\th_c\in (0,\pi)\cup (\pi,2\pi)$ and $\th_c/\pi\in \Q$. We set
$\sg_c=\th_c/2\pi\in(0,1)\cap(\Q\setminus\{\frac{1}{2}\})$. The closed
geodesic $c$ is called a degenerate saddle, if it satisfies
\be i(c)=1\quad {\it and}\quad k_0(c^{n_c}) = k_2(c^{n_c}) = 0, \lb{4.1}\ee
where $n_c$ and $k_j(c^{n_c})$ are defined in Definition 3.9.}

The following consequence of a result by H.-B. Rademacher in \cite{Rad2} will
be crucial in Section 10. We denote by $\alpha_c$ the mean index
$\hat{i}(c)\equiv\lim\limits_{m\rightarrow\infty} \frac{i(c^m)}{m}$ of a
closed geodesic $c$.

{\bf Theorem 5.2.} {\it Let $(S^2,F)$ be a Finsler $2$-sphere and
assume that there exists only one prime closed geodesic $c$
on $(S^2, F)$ and that $c$ is a rationally elliptic degenerate saddle. Then
there holds }
\be
\frac{n_c-1-\hat{k}_1(c^{n_c})}{n_c\aa_c} = 1. \lb{4.2}
\ee
The number $\hat{k}_1(c^{n_c})$ is defined in Definition 3.9.

{\bf Proof.} Theorem 7.9 in \cite{Rad2} treats compact, simply connected
Finsler manifolds $(M,F)$ with only finitely many prime closed geodesics and
provides a relation between invariants of these closed geodesics and a
topological invariant of $M$. A simple computation shows that this last
invariant equals $-1$ in the case $M=S^2$. If $c$ is the only prime closed
geodesic on $(S^2,F)$, this relation says
\be \frac{\beta_c}{\alpha_c} = -1,  \lb{4.3}\ee
where $\beta_c$ is the invariant of $c$ defined in \cite{Rad2}, Satz 7.3.
The proof of Theorem 5.2 consists in expressing $\beta_c$ by $n_c$
and $\hat{k}_1(c^{n_c})$ as follows
\be \beta_c=\frac{\hat{k}_1(c^{n_c}) + 1 - n_c}{n_c}.  \lb{4.4}\ee
Then (\ref{4.2}) is a direct consequence of (\ref{4.3}) and (\ref{4.4}).
The invariant $\beta_c$ is defined as follows. Set
$$  M_{m,j}(c) = \dim \overline{C}_j (E,c^m)  $$
with $\ol{C}_j (E,c^m)=H_j(\Lm(c^m)\cup S^1\cdot c^m/S^1,\Lm(c^m)/S^1)$,
cf. (\ref{2.8}). There exists a minimal even integer $k(c)>0$ such that for
all $m\ge 1, j\ge 0$:
$$ M_{m, j+i(c^m)} (c) = M_{m+k(c), j+i(c^{m+k(c)})}(c).  $$
Then
\be \beta_c
   = \frac{1}{k(c)}\sum\limits_{1\le m\le k(c), j\ge 0}(-1)^j M_{m,j}(c).
   \lb{4.5}\ee
Using Theorem 3.11, Proposition 3.12 and the index iteration formula (\ref{3.7})
we first compute the numbers $M_{m,j}(c)$ for a rationally elliptic degenerate
saddle $c$. According to (\ref{3.7}) we have
\bea
i(c^m) = 2[m\sg_c]+1, \quad \nu(c^m)=0, &\qquad&
                  {\rm if}\;m\not\in n\N, \lb{4.6}\\
i(c^m) = 2[m\sg_c]-1, \quad \nu(c^m)=2, &\qquad&
                  {\rm if}\;m\in n\N, \lb{4.7}
\eea
where $\sg_c=\th_c/2\pi\in(0,1)\cap\left(\Q\bs\{\frac{1}{2}\}\right)$
and $n=n_c\ge 3$ is the denominator of the reduced fraction $\sg_c$.

In particular, (\ref{4.6}) and (\ref{4.7}) imply
\be \epsilon(c^m)=(-1)^{i(c^m)-i(c)}=1  \lb{4.8}\ee
for all $m\in\N$. Using Proposition 3.12 we obtain
\be
M_{m,j}(c) = \left\{\matrix{
      1 \quad {\rm if}\;\;m\notin n\N \;\; {\rm and} \;\;j=i(c^m)  \cr
      0 \quad {\rm if}\;\;m\notin n\N \;\; {\rm and} \;\; j\not= i(c^m). \cr}\right.
      \lb{4.9}\ee

In the remaining cases we can use Theorem 3.11 and Proposition 3.12 to conclude
\be
M_{m,j}(c) = \left\{\matrix{
      \hat{k}_1 (c^n) \quad {\rm if}\;\;m\in n\N \;\; {\rm and} \;\;j=i(c^m)+1  \cr
   \quad    0 \quad \quad \;{\rm if}\;\;m\in n\N \;\; {\rm and} \;\; j\not= i(c^m)+1. \cr}\right.
      \lb{4.10}\ee

In particular, we can take $k(c)=2n$. Now $\beta_c$ can be computed from
(\ref{4.5}), and the result of this computation is equation (\ref{4.4}). \hb

Below, we will present a generalization of Theorem 5.2 which can
be used in the study of non-degenerate closed geodesics in Cases
CG-8 and CG-9 too. Following Theorem 7.9 of \cite{Rad2}, we have

{\bf Theorem 5.3.} {\it Let $S^2=(S^2,F)$ be a Finsler $2$-sphere
with only finitely many prime closed geodesics, each of which is
either non-degenerate, or a rationally elliptic degenerate saddle.
Denote non-degenerate prime closed geodesics on $S^2$ by $c_j$ for
$1\le j\le r$, and rationally elliptic degenerate saddle prime closed
geodesics on $S^2$ by $c_j$ for $r+1\le j\le r+a<+\infty$. We denote
the mean index of $c_j$ by $\aa_j=\hat{i}(c_j)$. Let
$\ga_j=\ga_{c_j}\in \{\pm 1/2, \pm 1\}$ such that
\be 2\ga_j = i(c_j^2)-i(c_j)\quad \mod\quad 2, \quad {\it and}
    \quad \ga_j(-1)^{i(c_j)}>0.  \lb{4.11}\ee
Then there holds
\be \sum_{j=1}^r\frac{\ga_j}{\aa_j}
 - \sum_{j=r+1}^{r+a}\frac{(n_j-1-\hat{k}_1(c^{n_j}_j))}{n_j\aa_j}
 = -1, \lb{4.12}\ee
where we set $n_j=n_{c_j}$ for $r+1\le j\le r+a$. }

{\bf Proof.} Note first that by (\ref{4.6}) and (\ref{4.7}), for a
rationally elliptic degenerate saddle closed geodesic $c$, we always
have $\ga_c=-1$.

Note that for a closed geodesic $c$, its mean index satisfies
either $\hat{i}(c)>0$ or $\hat{i}(c)=0$. When $\hat{i}(c)=0$, there holds
$i(c^m)=0$ for all $m\ge 1$.

Therefore for a non-degenerate elliptic or a rationally elliptic degenerate
saddle closed geodesic, because the rotation angle is always positive, by
the iteration formula (\ref{3.7}) or (\ref{3.8}), its mean index is
positive.

For any hyperbolic closed geodesic $c$, if $\hat{i}(c)=0$, then $i(c^m)= 0$
for all $m\ge 1$. Note that there holds $E(c^m)>E(c^{m-1})$ for all $m\ge 2$.
Because $c^m$ is non-degenerate in $\Lm$, each $c^m$ for $m\ge 1$ must be a
strict local minimum of $E$ in $\Lm$. For every strict local minimum $c^m$
with $m\ge 2$, we have $E(c^{m-1})<E(c^m)$. Thus a mountain pass argument
yields a closed geodesic $d_m$ with a non-trivial $1$-dimensional local
homological critical module and $E(d_m)>E(c^m)$. This argument yields
infinitely many such $d_m$s. Note that here we do not claim that these
$d_m$s yield infinitely many prime closed geodesics. But if the total number
of prime closed geodesics is finite, indices of iterations of each one
should either increase to infinity or be identically equal to zero
respectively. Then their iterates can not possess enough non-trivial
$1$-dimensional local homological critical modules by Proposition 3.8.
This contradiction proves $\hat{i}(c)>0$.

Therefore the denominators on the left hand side of (\ref{4.12}) are
all non-zero, and then (\ref{4.12}) holds.

Now using our Theorem 5.2 together with the proof of Theorem 3.1 of
\cite{Rad1} and \cite{Rad2}, we get Theorem 5.3. \hb

\setcounter{equation}{0}
\section{Two theorems by N. Hingston}

In certain cases one can use the existence of a degenerate closed geodesic
$c$ to prove the existence of infinitely many closed geodesics. The first
instance of this phenomenon was discovered in \cite{Ban1} in the case of
vanishing mean index $\hat{i}(c)=0$, see also \cite{Ban3}. The method
was considerably advanced by N. Hingston in \cite{Hin1} and \cite{Hin2} who
was able to treat cases where $\hat{i}(c)>0$. The proofs of these results
are variational and do not require symmetry of the metric. Hence the results
apply to general Finsler metrics. The following statement combines
\cite{Hin1}, Proposition 1, and \cite{Hin2}, Theorem, for the case of a
Finsler 2-sphere $(S^2, F)$.

{\bf Theorem 6.1.} {\it Let $c$ be a closed geodesic on $(S^2,F)$ that satisfies
(Iso). Assume that either }
\begin{itemize}
\item[(1)]{\it $k_0(c)>0$ and $i(c^m)=m (i(c)+1)-1$, $\nu(c^m)=\nu(c)$ for
all $m\in\N$, or}
\end{itemize}
\begin{itemize}
\item[(2)]{\it $k_{\nu(c)}(c)>0$ and $i(c^m)+\nu(c^m) =m(i(c)+\nu(c)-1)+1$,
$\nu(c^m)=\nu(c)$ for all $m\in\N$.}
\end{itemize}
{\it Then there exist infinitely many prime closed geodesics on $(S^2,F)$.}

To see that \cite{Hin1}, Proposition 1, and \cite{Hin2}, Theorem, imply
Theorem 6.1, note that Lemma 3.3, (\ref{2.25}) and (\ref{2.26})
imply $C_{i(c)}(E,c)\not= 0$ in case (1), and $C_{i(c)+\nu(c)}(E,c)\not= 0$
in case (2). So, in case (1) the hypotheses of \cite{Hin2}, Theorem, are
satisfied, while in case (2) the hypotheses of \cite{Hin1}, Proposition 1,
hold.

By Lemma 3.10 the conditions $k_0(c)>0$ or $k_{\nu(c)}(c)>0$, respectively,
are equivalent to the fact that $c$ is a strict local minimum or a strict local
maximum of $E$ restricted to a local characteristic manifold at $c$.

N. Hingston imposes the seemingly weaker assumptions $i(c^m)\ge m(i(c)+1)-1$
in case $(1)$, and $i(c^m)+\nu(c^m)\le m(i(c)+\nu(c)-1)+1$ in case (2).
The estimates in Theorems 10.1.2 and 10.1.3 of \cite{Lon3}, originally proved
in \cite{LLo1}, imply that these inequalities are in fact equalities.

\setcounter{equation}{0}
\section{Homology of $(\Lm S^2,\Lm^0S^2)$ and first consequences }

In \cite{Zil1} W. Ziller computed the $\Z$-homology of the free loop space
$\Lambda$ of compact rank 1 symmetric spaces (with the exception of $\R P^n$).
For $\Lambda S^2$ the table on p. 21 of \cite{Zil1}, taken literally, does not
give the correct result. However, the correct result follows easily from
\cite{Zil1}, Theorem 8, and it is explicitly stated in \cite{Zil2}, p. 148.
Specialized to the coefficient ring $\Q$ one has
\be  H_k(\Lambda S^2)=\Q  \lb{6.1}\ee
for all $k\in\N\cup\{0\}$.

Next we solve the simple exercise to compute $H_k(\Lm S^2, \Lm^0S^2)$ from
(\ref{6.1}). If $i:\Lambda^0 S^2\rightarrow \Lambda S^2$ denotes inclusion
and $\ev:\Lambda S^2\rightarrow S^2$ denotes the evaluation map
$\ev(\gamma)=\gamma(1)$, then $\ev\circ i:\Lambda^0 S^2\rightarrow S^2$
is a diffeomorphism. This implies that
$i_* : H_* (\Lm^0 S^2) \rightarrow H_*(\Lm S^2)$ is one-to-one. Hence the long
exact homology sequence of the pair $(\Lambda S^2, \Lambda^0 S^2)$ together
with (\ref{6.1}) show that
\be
H_k(\Lm S^2,\Lm^0S^2) = \left\{\matrix{
     0, & \quad {\it if}\;\;k= 0,\;\;{\it or}\;\;k=2, \cr
     \Q, & \quad {\it if}\;\;k=1, \;\;{\it or}\;\;k\ge 3. \cr}\right.
\lb{6.2}\ee

From now on in the rest of this paper, we write simply $\Lm=\Lm S^2$ and
$\Lm^a=\Lm^a S^2$ for $a\in\R$. In the following three sections we will prove
Theorem 1.1 by contradiction. So we will assume the condition
\newline
{\bf (F) There exists only one prime closed geodesic $c$ on the given Finsler
2-sphere $(S^2, F)$. }

We mention some simple consequences of $(F)$, (\ref{6.2}) and the Morse
inequalities (\ref{2.33}). Using Proposition 3.8 and (\ref{6.2}) and the fact
that $\nu(c^m)\le 2$ for all $m\in\N$, we see that the sequence $i(c^m)$ is
unbounded. By the index iteration formulae (\ref{3.1})-(\ref{3.9}) this
implies $i(c^m)\ge i(c)>0$ for all $m\in\N$. Using Proposition 3.8 again, we
conclude that the Morse type number $M_0$ satisfies
\be  M_0=0.  \lb{6.3}\ee
Moreover, by (\ref{2.32}) and (\ref{6.2}) we have
$$  M_1\ge b_1(\Lambda, \Lambda^0)=1,  $$
and hence $i(c)=1$. By (\ref{3.3}) this implies:
\be \mbox{The only prime closed geodesic $c$ cannot be of type CG-3.}
       \lb{6.4}\ee

Moreover, the integer $p$ in the iteration formulae (\ref{3.1}), (\ref{3.2})
and (\ref{3.4})-(\ref{3.9}) satisfies:
\bea
&&\mbox{If $c$ is of one of the types CG-1, CG-2 or CG-9, then $p=1$;}
       \lb{6.5}\\
&&\mbox{If $c$ is of one of the types CG-4, CG-5, CG-6, CG-7 or CG-8, then $p=0$.}
       \lb{6.6}
\eea

\setcounter{equation}{0}
\section{Cases with eigenvalue $1$}

We recall that we use homology with rational coefficients. We shall use
the results from Section 3 to compute local critical modules. We recall
the numbers $k_j(c)$ and $\hat{k}_j(c)$ defined in Definition 3.9.

\subsection{Case CG-$k$ with $k=1$ or $k=2$}

Note that (\ref{3.1}) and (\ref{3.2}) imply that $i(c^m)-i(c)$ is even for
every $m\in\N$. Moreover, (\ref{6.5}) implies that $p=1$ in formulae
(\ref{3.1}) and (\ref{3.2}). Thus, if $c$ is of type CG-$k$ with
$k\in\{1,2\}$, then (\ref{3.1}) and (\ref{3.2}) become
\be
i(c^m)=2m-1, \quad \nu(c^m)=k
         \qquad {\rm for \,\, all}\;\; m\ge 1. \lb{7.1}
\ee
For the closed geodesic $c$ itself we obtain from Proposition 3.8
\be C_1(E, S^1\cdot c) = H_0 (N^-_c\cup\{c\}, N^-_c)^{\Z_m}=\Q^{k_0(c)},
  \lb{7.2}\ee
since $\hat{k}_0(c)=k_0(c)$, cf. (\ref{2.27}).

For the iterates $c^m$ with $m\ge 2$, we have $i(c^m)\ge 3$, so that
Proposition 3.8 implies
\be  C_1(E, S^1\cdot c^m)=0 \qquad {\rm for }\;\; m\ge 2. \lb{7.3}\ee

Using the Morse inequality (\ref{2.32}) and the fact that
$b_1=b_1(\Lambda, \Lambda^0)=1$ by (\ref{6.2}), and (\ref{7.2}) and
(\ref{7.3}), we obtain
\be  k_0(c)=M_1\ge b_1=1. \lb{7.4}\ee

Now (\ref{7.1}) and (\ref{7.4}) imply that the hypothesis (1) of
N. Hingston's Theorem 6.1 is satisfied. Hence, in contradiction to our
assumption (F), there exist infinitely many prime closed geodesics on
$(S^2, F)$.

\subsection{Case CG-3}

According to (\ref{6.4}) this case cannot occur.

\setcounter{equation}{0}
\section{Cases with eigenvalue $-1$}

\subsection{Case CG-$k$ with $k=4$ or $k=5$}

In these two cases, $i(c^m)-i(c)$ is even for every $m\in\N$ by
(\ref{3.4}) or (\ref{3.5}). According to (\ref{6.6}), we have $p=0$ in the
formulae (\ref{3.4}) and (\ref{3.5}). Then (\ref{3.4}) or (\ref{3.5})
become
\be i(c^m)=m-\frac{1+(-1)^m}{2}, \quad
   \nu(c^m)= \frac{(1+(-1)^m)[(k-1)/2]}{2},
   \qquad \mbox{ for all }\quad m\ge 1. \lb{8.1}\ee

We will now compute the Morse type numbers $M_1, M_2$ and $M_3$.

Since $i(c)=1$ and $\nu(c)=0$ we obtain from Proposition 3.8
\bea
C_1(E,S^1\cdot c) &=& C_2(E,S^1\cdot c) =\Q, \lb{8.2}\\
C_3(E,S^1\cdot c) &=& 0.  \lb{8.3}
\eea

For $c^2$, we have $i(c^2)=1=i(c)$, $\nu(c^2)=[(k-1)/2]$, and $i(c^2)-i(c)=0$.
Then by Proposition 3.8, we obtain
\bea
C_1(E,S^1\cdot c^2) &=& \Q^{\hat{k}_0(c^2)}, \lb{8.4}\\
C_2(E,S^1\cdot c^2) &=& \Q^{\hat{k}_1(c^2)+\hat{k}_0(c^2)}, \lb{8.5}\\
C_3(E,S^1\cdot c^2) &=& \Q^{\hat{k}_2(c^2)+\hat{k}_1(c^2)}. \lb{8.6}
\eea

For $c^3$, we have $i(c^3)=3$, $\nu(c^3)=0$, and $i(c^3)-i(c)=2$. Thus
by Proposition 3.8, we have
\bea
C_1(E,S^1\cdot c^3) &=& C_2(E,S^1\cdot c^3) = 0, \lb{8.7}\\
C_3(E,S^1\cdot c^3) &=& \Q. \lb{8.8}
\eea

For $c^4$, we have $i(c^4)=3$, $\nu(c^4)=1$ in Case CG-4 and $\nu(c^4)=2$
in Case CG-5.
Because $i(c^4)-i(c^2)=2$, by Proposition 3.8,
in both cases CG-4 and CG-5 we have
\bea
C_1(E,S^1\cdot c^4) &=& C_2(E,S^1\cdot c^4) = 0, \lb{8.9}\\
C_3(E,S^1\cdot c^4) &=& \Q^{\hat{k}_0(c^4)}.  \lb{8.10}
\eea

For $c^m$ with $m\ge 5$, we have $i(c^m)\ge 4$ and hence Proposition 3.8 implies
\be C_q(E,S^1\cdot c^m) = 0, \qquad {\rm for\;all}\; q\le 3, \;\;m\ge 5. \lb{8.11}
\ee

Thus we obtain $M_1=1+\hat{k}_0(c^2)$, $M_2=1+\hat{k}_0(c^2)+\hat{k}_1(c^2)$,
and $M_3=1+\hat{k}_0(c^4)+\hat{k}_1(c^2)+\hat{k}_2(c^2)$.
Then by (\ref{2.33}), (\ref{6.2}), and (\ref{6.3}), we obtain
\be  1+\hat{k}_0(c^4)+\hat{k}_2(c^2) = M_3-M_2+M_1
       \ge b_3-b_2+b_1 = 2.   \lb{8.12}\ee

Therefore, at least one of $\hat{k}_0(c^4)$ and $\hat{k}_2(c^2)$ must be
positive. First suppose that $\hat{k}_0(c^4)$ is positive. Then, by
(\ref{2.27}) and Theorem 3.11, we have
$$  k_0(c^2) = \hat{k}_0(c^2) = \hat{k}_0(c^4)>0.  $$
Now consider the closed geodesic $d=c^2$. Then we have $k_0(d)>0$,
$i(d)=1$, and
\be i(d^m)=i(c^{2m})=2m-1=m(i(d)+1)-1, \quad \nu(d^m)=\nu(d), \quad
                {\rm for\;all}\; m\ge 1. \lb{8.13}\ee
Therefore, by Theorem 6.1, there exist infinitely many prime closed geodesics
on $(S^2, F)$.

Finally suppose that $\hat{k}_2(c^2)$ is positive. This can only happen in
case CG-$5$ when $\nu(c^2)=2$. Considering $d=c^2$ again we have
$k_{\nu(d)}(d)>0$, $i(d)=1$, $\nu(d)=2$, and for all $m\ge 1$:
\bea
&& i(d^m) = i(c^{2m}) = 2m-1 = m(i(d)+1)-1, \quad
     \nu(d^m)=\nu(d), \lb{8.14}\\
&& i(d^m) + \nu(d^m) = 2m+1 = m(i(d)+\nu(d)-1)+1.   \lb{8.15}
\eea
Again, by Theorem 6.1, we obtain infinitely many prime closed geodesics
on $(S^2, F)$.

\subsection{Case CG-6}

According to (\ref{6.6}) we have $p=0$ in formula (\ref{3.6}), so that
(\ref{3.6}) becomes
\be
i(c^m)=m, \quad \nu(c^m)= \frac{1+(-1)^m}{2},
     \qquad {\rm for\;all}\; m\ge 1.  \lb{8.16}\ee
Note that
\be \epsilon(c^m)=(-1)^{i(c^m)-i(c)}=\left\{\matrix{
    -1, &\quad {\rm if}\quad m\quad {\rm is}\quad {\rm even}, \cr
    +1, &\quad {\rm if}\quad m\quad {\rm is}\quad {\rm odd}. \cr}
\right.\lb{8.17}\ee
Next we compute the Morse type numbers $M_1$, $M_2$ and $M_3$.

Since $i(c)=1$ and $\nu(c)=0$, Proposition 3.8 implies
\be \matrix{
C_1(E, S^1\cdot c) = C_2(E, S^1\cdot c)=\Q,    \cr
C_3(E, S^1\cdot c)= 0.                   \cr}\lb{8.18}\ee
From (\ref{8.16}) we have $i(c^2)=2$, $\ep(c^2)=-1$ and $\nu(c^2)=1$. So
Proposition 3.8 and (\ref{8.17}) imply
$$  C_q(E, S^1\cdot c^2) =
  H_{q-2}\left(N^-_{c^2}\cup\{c^2\}, N^-_{c^2}\right)^{\Z_2, -1}
   \oplus H_{q-3}\left(N^-_{c^2}\cup\{c^2\}, N^-_{c^2}\right)^{\Z_2, -1}, $$
where $N_{c^2}$ denotes a local characteristic manifold at $c^2$,
$N^-_{c^2}=N_{c^2}\cap\Lambda(c^2)$. We will show that
$H_*\left(N^-_{c^2}\cup\{c^2\}, N^-_{c^2}\right)^{\Z_2, -1}=0$, and this
will prove
\be  C_*(E,S^1\cdot c^2) = 0. \lb{8.19}\ee

First note that $H_0\left(N^-_{c^2}\cup\{c^2\}, N^-_{c^2}\right)^{\Z_2, -1}=0$,
since the $\Z_2$-action fixes $c^2$. From the assumption
$H_1\left(N^-_{c^2}\cup\{c^2\},N^-_{c^2}\right)^{\Z_2, -1}\not= 0$, we
conclude that $k_1(c^2)=\dim H_1(N^-_{c^2}\cup\{c^2\}, N^-_{c^2})>0$.
In this case the hypotheses (2) of Theorem 6.1 hold for the closed geodesic
$d=c^2$: We have $\nu(d)=1, k_{\nu(d)}(d)>0$ and, by (\ref{8.16}),
$$ i(d^m)+\nu(d^m)=2m+1=m(i(d)+\nu(d)-1) +1, \quad \nu(d^m)=\nu(d) $$
for all $m\ge 1$. Thus we obtain infinitely  many prime closed geodesics,
in contradiction to our assumption $(F)$.
Hence we have $H_1\left(N^-_{c^2}\cup\{c^2\}, N^-_{c^2}\right)^{\Z_2, -1}= 0$.
Since $\dim N_{c^2}=\nu(c^2)=1$, this completes the proof that
$H_*(\left(N^-_{c^2}\cup\{c^2\}, N^-_{c^2}\right)^{\Z_2, -1}= 0$.

Since $i(c^3)=3$ and $\nu(c^3)=0$, Proposition 3.8 implies
\be\matrix{
C_1(E, S^1\cdot c^3) = C_2(E, S^1\cdot c^3)=0  \cr
C_3(E, S^1\cdot c^3)= \Q.                   \cr}
\lb{8.20}\ee
Finally, if $m\ge 4$ then $i(c^m)\ge 4$, and hence by Proposition 3.8 we
obtain
\be C_q(E, S^1\cdot c^m)=0   \qquad {\rm for }\;\; q\in\{1,2,3\}.
\lb{8.21}\ee

Now (\ref{8.18})-(\ref{8.21}) imply $M_1=M_2=M_3=1$. From the Morse
inequalities (\ref{2.33}) and from (\ref{6.2}), (\ref{6.3}),
we obtain the contradiction
$$  1=M_3-M_2+M_1\ge b_3-b_2+b_1=2.  $$
Hence only one closed geodesic of type CG-6 cannot generate all the
homology of $(\Lm, \Lm^0)$.

\setcounter{equation}{0}
\section{Case CG-7 of a rationally elliptic closed geodesic}

In this section we will derive a contradiction from the assumption (F)
that a Finsler sphere $(S^2, F)$ has only one prime closed geodesic $c$
if this $c$ is of type CG-7, i.e., if the linearized Poincar\'e map of $c$ is
conjugate to a rotation by an angle $\th\in(0,\pi)\cup(\pi,2\pi)$
with $\th/\pi\in\Q$.

Our arguments use first N. Hingston's results, Theorem 6.1, to reduce the
problem to the subcase of a rationally elliptic degenerate saddle, cf.~Definition 5.1.
Then using Theorem 5.2 or Theorem 5.3 which are based on
H.-B. Rademacher's work and the index iteration
formula (\ref{3.7}), we further restrict the rotation angle $\th$. These
results allow us to show that $c$ and its iterates generate a
surplus in local one-dimensional homology. A careful analysis of the
situation then shows that the local 2-dimensional homology generated by
the iterates of $c$ cannot destroy this surplus in one-dimensional
homology. This is the final contradiction. More precisely, this
contradiction is reached in the following three steps.

\smallskip

{\bf Step 1.} {\it General information on the closed geodesic $c$}

We first mention some consequences of our assumptions on $c$.

We set
\be \sg=\th/2\pi. \lb{9.1}\ee
Then $\sg\in(0,1)\cap(\Q\bs\{1/2\})$. From (\ref{6.6}) we know that
the integer $p$ in (\ref{3.7}) equals zero. Hence (\ref{3.7}) becomes
\bea
i(c^m)=2[m\sg]+1,\quad \nu(c^m)=0,
          &\qquad& {\rm if}\;\;m\sg\not\in \N,  \lb{9.2}\\
i(c^m)=2[m\sg]-1,\quad \nu(c^m)=2,
          &\qquad& {\rm if}\;\;m\sg \in \N. \lb{9.3}
\eea

In particular, we have $i(c)=1$ and the mean index $\aa\equiv\hat{i}(c)$
satisfies $\aa=2\sg$. Moreover, with $n\equiv n_c\in\N$ given by
Definition 3.9, we have $n\ge 3$, and $k\equiv n\sg$ satisfies
$k\in\{1,\ldots,n-1\}$ and is relatively prime to $n$. Then $d=c^n$ is a
degenerate closed geodesic satisfying $i(d)=2k-1$, $\nu(d)=2$,
$\epsilon(c^n)=(-1)^{i(c^n)-i(c)}=1$, and
\bea
&& i(d^m)=2km-1 = m(i(d)+1)-1, \qquad \nu(d^m)=2,   \lb{9.4}\\
&& i(d^m)+\nu(d^m)=2km+1=m(i(d)+\nu(d)-1)+1, \lb{9.5}
\eea
for all $m\in\N$. Hence, the closed geodesic $d$ satisfies part of the
assumptions of N. Hingston's Theorem 6.1. Since this result promises the
existence of infinitely many closed geodesics on $(S^2,F)$, its
additional assumptions $k_0(d)>0$ or $k_2(d)>0$ are both not true, i.e.,
$c$ is a rationally elliptic degenerate saddle in the sense of Definition
5.1. Hence our Theorem 5.2 implies the following crucial identity, which
restricts the rotation angle $\th$:
\be  n-1 -\hat{k}_1(c^n) = n \aa,   \lb{9.6}\ee
where $\hat{k}_1(c^n)$ is given by Definition 3.9.

Since $\aa=2\sg$ and $\sg=k/n$ with $k\ge 1$, we obtain
\be 0\le\hat{k}_1(c^n)=n-1-2k\le n-3,   \lb{9.7}\ee
and, in particular, $2k+1\le n$. Hence we have
\be 2\sg <1, \lb{9.8}\ee
i.e., $\th<\pi$. We set
\be  \tau=\max\{m\in\N \mid m\sg<1\}. \lb{9.9}\ee
Because of (\ref{9.8}) and $\sg=k/n$ we have
\be 2\le\tau\le n-1. \lb{9.10}\ee

\smallskip

{\bf Step 2.} {\it Vanishing connecting homomorphism and additive
homologies among level sets}

Now we will study the one-dimensional homology generated by the closed
geodesics $c$, $c^2$, $\ldots$, $c^{\tau}$. By (\ref{9.2}) and (\ref{9.9})
all of them are non-degenerate and of index one. We set $\ka_0=0$ and
$$ \ka_m=E(c^m), \qquad {\rm for\;all}\; m\in\N. $$
There holds
\bea
&&  0=\ka_0<\ka_1<\cdots \ka_m<\ka_{m+1}<\cdots, \lb{9.11}\\
&&  \ka_m\to +\infty, \qquad {\rm as}\quad m\to +\infty. \lb{9.12}
\eea
Note that
\be H_0(\Lm^{\ka_m},\Lm^0)=0, \qquad {\rm for\;all}\; m\in\N, \lb{9.13}\ee
holds, since there are no closed geodesics of index zero. We recall
that, for all $q\in\N\cup\{0\}$ and $m\in\N$, we have
$H_q(\Lm^{\ka_m},\Lm^{\ka_{m-1}})=C_q(E,S^1\cdot c^m)$, cf.
Proposition 3.13.

In \cite{BoS1} of 1958, R. Bott and H. Samelson established the
additivity of homologies of level sets for pointed loop spaces of
compact globally symmetric spaces. In \cite{Zil1} of 1977, W. Ziller
established this additivity for free loop spaces of such spaces.
In general $(S^2,F)$ is not a globally symmetric space. Our following
result establishes also such an additivity of the homologies of level
sets of the energy functional on $(\Lm^{\ka_{\tau}},\Lm^0)$ for Case
CG-7 by a rather different method.

{\bf Proposition 10.1.} {\it Under the assumption (F), let $c$ be a
prime closed geodesic of type CG-7. Then }
$$ H_1(\Lm^{\ka_{\tau}}, \Lm^0) = \Q^{\tau}. $$

{\bf Proof.} In Lemma 10.2 below we will show that for every $2\le m\le\tau$
the connecting homomorphism
$$ \pt_2:H_2(\Lm^{\ka_m},\Lm^{\ka_{m-1}})
    \rightarrow H_1(\Lm^{\ka_{m-1}},\Lm^0)  $$
of the exact homology sequence of the triple
$(\Lm^{\ka_m},\Lm^{\ka_{m-1}},\Lm^0)$ is zero. Using this and the fact
that $H_0(\Lm^{\ka_{m-1}},\Lm^0)=0$, cf. (\ref{9.13}), this exact
sequence splits and implies that
$$ H_1(\Lm^{\ka_m},\Lm^0) =
     H_1(\Lm^{\ka_{m-1}},\Lm^0)\oplus H_1(\Lm^{\ka_m},\Lm^{\ka_{m-1}}). $$
Since $i(c^m)-i(c)=0$, (i) of Proposition 3.8 implies that
$H_1(\Lm^{\ka_m},\Lm^{\ka_{m-1}}) = \Q$ for $1\le m\le\tau$. Hence
our claim follows by induction. \hb

{\bf Lemma 10.2.} {\it Under the assumption of Proposition 10.1, for
$2\le m\le\tau$ the connecting homomorphism
$$ \pt_2: H_2(\Lm^{\ka_m},\Lm^{\ka_{m-1}})
        \rightarrow H_1(\Lm^{\ka_{m-1}},\Lm^0)  $$
of the triple $(\Lm^{\ka_m},\Lm^{\ka_{m-1}},\Lm^0)$ is zero.}

{\bf Proof.} Recall that
\be i(c^m)=1, \quad \nu(c^m)=0, \qquad {\rm for}\; 1\le m\le\tau.
           \lb{9.14}\ee
So the local negative disks $U_{c^m}$ are one-dimensional, and we can
assume that the m'th iteration map $\phi_m$ maps $U_c$ onto $U_{c^m}$.
This implies that the $\Z_m$-action on $U_{c^m}$ is trivial. Using
Proposition 3.8 and (\ref{9.14}) we see that
\be H_2(\Lm^{\ka_m},\Lm^{\ka_{m-1}}) = C_2(E,S^1\cdot c^m)
    = \Q, \qquad {\rm for\;} 1\le m\le\tau. \lb{9.15}\ee
Note that a representative for a generator of
$H_2(\Lm^{\ka_m},\Lm^{\ka_{m-1}}) = C_2(E,S^1\cdot c^m)$ can be
constructed as follows.

(i) Let $f:([-1,1],\{-1,1\}) \rightarrow (\Lm^{\ka_1},\Lm^0)$ be a
continuous map such that $f(0)=c$ and satisfying $0<E(f(t))<\ka_1$ for
$t\in (-1,1)\bs\{0\}$, and such that, for some $\epsilon>0$,
$f|_{[-\epsilon,\epsilon]}$ represents a generator of $H_1(U_c, U_c^-)$.
Such an $f$ exists by (\ref{9.14}). Then, by Theorem 3.11, for every
$m\in\{1,2,\ldots,\tau\}$, the curve
$f^m:=\phi_m\circ f:[0,1]\rightarrow\Lambda^{\kappa_m}$ has the property that
$f^m|_{[-\epsilon,\epsilon]}$ generates $H_1(U_{c^m}, U^-_{c^m})$, cf.
Figure 10.1.

(ii) Define $F_m:(S^1\times[-1,1],S^1\times\{-1,1\}) \rightarrow
(\Lm^{\ka_m},\Lm^{\ka_{m-1}})$ by $F_m(\th,t)=\th\cdot f^m(t)$, cf.
Figure 10.1. If $0\not= h\in H_2(S^1 \times [-1,1], S^1 \times \{-1,1\})$
is the standard generator, then
$0\not= (F_m)_*(h)\in H_2(\Lm^{\ka_m}, \Lm^{\ka_{m-1}})$,
cf. the proof of Proposition 3.2. Hence $(F_m)_*(h)$ generates
$H_2(\Lm^{\ka_m},\Lm^{\ka_{m-1}})$.

\begin{figure}
\begin{center}
\resizebox{9.5cm}{9.5cm}{\includegraphics*[1cm,1cm][17cm,17cm]{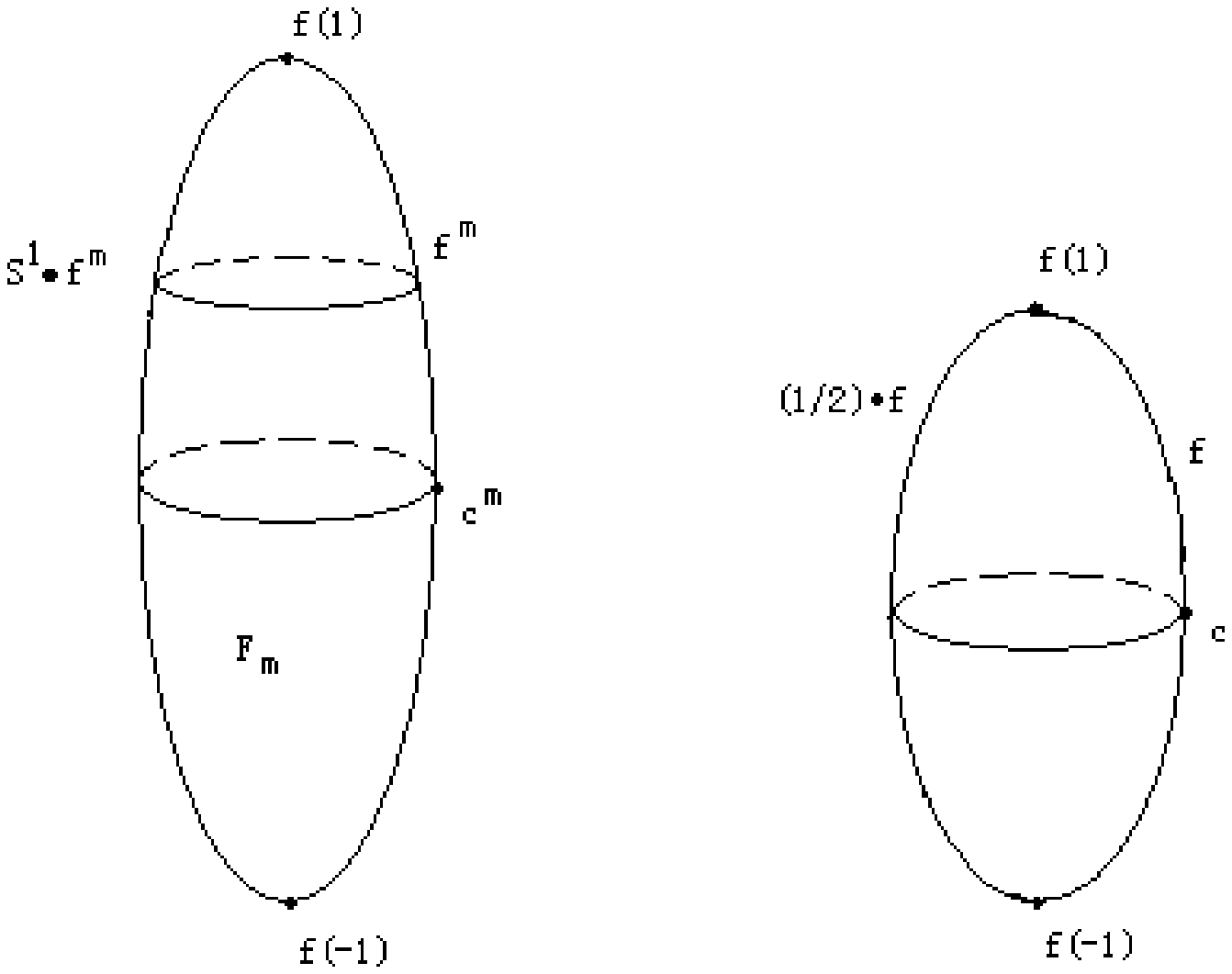}}
\caption{The map $F_m$}
\end{center}
\end{figure}
\vspace{2mm}

By our above construction, $\pt_2((F_m)_* h)=(F_m)_*(\pt_2 h)$
is the difference of the homology classes of two trivial $S^1$-orbits
\bea
\pt_2((F_m)_* h)
&=& [S^1\cdot f^m(1)] - [S^1\cdot f^m(-1)] \nn\\
&=& [f(1)] - [f(-1)]   \nn\\
&=& 0 \qquad {\rm in}\quad H_1(\Lm^{\ka_{m-1}},\Lm^0). \quad \mbox{\hb}
\nn\eea

\smallskip

{\bf Step 3.} {\it Chasing exact sequences}

Now we treat the case where $\tau<n-1$, i.e., $k=n\sg>(\tau+1)\sg>1$.
By (\ref{9.2}) and (\ref{9.3}) this entails $i(c^m)\ge 3$ for all $m>\tau$.
Hence we have
$$ H_2(\Lm,\Lm^{\ka_{\tau}})=0. $$
Now the exact sequence
$$ H_2(\Lm,\Lm^{\ka_{\tau}})\rightarrow H_1(\Lm^{\ka_{\tau}},\Lm^0)
      \rightarrow H_1(\Lm,\Lm^0) $$
implies that $\dim H_1(\Lm^{\ka_{\tau}},\Lm^0) \le \dim H_1(\Lm,\Lm^0)$,
where $\dim H_1(\Lm,\Lm^0)=1$ by (\ref{6.2}). This contradicts
Proposition 10.1 and the fact that $\tau\ge 2$, cf. (\ref{9.10}).

Finally, we consider the case that $\tau=n-1$, i.e., $k=n\sg=(\tau+1)\sg=1$.
Then we have $c^{\tau+1}=c^n=d$, $i(d)=1$, $\nu(d)=2$ and $k_0(d)=k_2(d)=0$,
since $c$ is a rationally elliptic degenerate saddle. Together with
Proposition 3.8, we obtain
\be \dim H_2(\Lm^{\ka_{\tau+1}},\Lm^{\ka_{\tau}}) = \hat{k}_1(d). \lb{9.16}\ee
Since $i(c^m)\ge 3$ for $m>\tau+1$, we have $H_2(\Lm,\Lm^{\ka_{\tau+1}})=0$.
Hence the exact sequence
$$ H_2(\Lm^{\ka_{\tau+1}},\Lm^{\ka_{\tau}}) \rightarrow
   H_2(\Lm,\Lm^{\ka_{\tau}}) \rightarrow H_2(\Lm,\Lm^{\ka_{\tau+1}}) $$
implies that
$\dim H_2(\Lm,\Lm^{\ka_{\tau}})
       \le \dim H_2(\Lm^{\ka_{\tau+1}},\Lm^{\ka_{\tau}})$.
Using (\ref{9.16}) we see that
\be \dim H_2(\Lm,\Lm^{\ka_{\tau}}) \le \hat{k}_1(d). \lb{9.17}\ee
The exact sequence
$$ H_2(\Lm,\Lm^{\kappa_{\tau}}) \rightarrow H_1(\Lm^{\ka_{\tau}},\Lm^0)
     \rightarrow H_1(\Lm,\Lm^0) $$
and $H_1(\Lm,\Lm^0)=\Q$, cf. (\ref{6.2}) imply
$$ \dim H_1(\Lm^{\ka_{\tau}},\Lm^0)\le \dim H_2(\Lm,\Lm^{\ka_{\tau}})+1. $$
Using Proposition 10.1 and (\ref{9.17}) we obtain
\be  \tau\le\hat{k}_1(d)+1. \lb{9.18}\ee
Since $\tau=n-1$, this contradicts (\ref{9.7}). \hb

Therefore we have proved that, under the assumption (F), the only prime closed
geodesic on $S^2$ cannot be of the class CG-7 with $k_0(c^n)=k_2(c^n)=0$.

{\bf Remark 10.3.} Note that, in Katok's example, for the two closed
geodesics $c_1$ and $c_2$ there holds
\be i(c_1)=1, \quad i(c_1^{m+1})\ge 3, \quad i(c_2^m)\ge 3, \quad
   {\rm for\;all}\; m\ge 1. \lb{9.19}\ee
Therefore in (\ref{9.9}) we have $\tau=1$
and $H_2(\Lm,\Lm^{\ka_1})=H_1(\Lm,\Lm^{\ka_1})=0$. Thus the long
exact homology sequence of the triple $(\Lm,\Lm^{\ka_{\tau}},\Lm^0)$
becomes
\bea
&& H_2(\Lm,\Lm^0)
  \mapright{} H_2(\Lm,\Lm^{\ka_{\tau}})
  \mapright{} H_1(\Lm^{\ka_{\tau}},\Lm^0)
  \mapright{} H_1(\Lm,\Lm^0)
  \mapright{} H_1(\Lm,\Lm^{\ka_{\tau}}) \nn\\
&& \qquad \parallel \qquad\qquad\qquad  \parallel \qquad\qquad\qquad
  \parallel \qquad\qquad\quad \parallel \qquad\qquad\quad \parallel \nn\\
&& \qquad 0 \qquad\qquad\quad\quad\, 0
  \qquad\qquad\quad\;\, \Q \qquad\qquad\;\;\, \Q \qquad\qquad\quad 0.
\lb{9.20}\eea
This is very different from the case of only one rationally
elliptic degenerate saddle closed geodesic.

\smallskip

Now finally we can give

{\bf Proof of Theorem 1.1.} Note that, by both Theorem 3.1 and Example 4.1
of \cite{Rad1}, it is impossible that the only prime closed
geodesic $c$ on $S^2$ in the assumption (F) in the Section 7
is of type CG-8 or CG-9 in Section 4. Here, for the reader's convenience,
we briefly indicate how to exclude these two cases. Note that by
$M_1\ge b_1=1$, we must have $i(c)=1$ in both cases. In Case CG-8, by
Theorem 5.3 we get $\aa_c=1$ which contradicts that $\aa_c$ should be
irrational. In Case CG-9, we get $i(c^m)=m$ and $\nu(c^m)=0$ for all
$m\ge 1$. Thus similarly to our study in Section 9.2, Proposition 3.8 and
a direct computation show $M_1=M_2=M_3=1$. This contradicts the Morse
inequality (\ref{2.33}), since $1=M_3-M_2+M_1\ge b_3-b_2+b_1=2$, cf.
(\ref{6.2}) and (\ref{6.3}).

The preceding Sections 7-10 show that under the assumption $(F)$ the only
prime closed geodesic $c$ cannot be of classes CG-1 to CG-7 either.
Therefore the proof of Theorem 1.1 is complete. \hb

\smallskip

{\bf Acknowledgements.} For his visits to Freiburg from February 9th to March
8th and July 3rd to 5th, 2004, Y. Long sincerely thanks the Mathematics
Institute of the University of Freiburg for its hospitality and the DFG for
its financial support. The authors sincerely thank Professor Jean-Pierre
Bourguignon for his careful reading of the draft and valuable comments on it.
They thank sincerely also Dr. Duanzhi Zhang and Dr. Wei Wang for helpful
discussions on the first draft of this paper, and Professor Kai Cieliebak
for drawing their attention to the number of closed geodesics mentioned in
Anosov's paper \cite{Ano1}.


\bibliographystyle{abbrv}

\bigskip

\noindent Victor Bangert,

Mathematisches Institut, Abteilung f\"ur Mathematik,
Albert-Ludwigs-Universit\"at, D-79104 \\ Freiburg im
Breisgau, Germany. \quad
E-mail: bangert@mathematik.uni-freiburg.de

\bigskip

\noindent Yiming Long,

Chern Institute of Mathematics and LPMC, Nankai University,
Tianjin 300071, the People's Republic of China. \quad
E-mail: longym@nankai.edu.cn

\bigskip

First version: May 2, 2006.

Revised version: April 16, 2007.

\end{document}